\documentclass{article}

\makeatletter
\@addtoreset{figure}{section}
\def\thefigure{\thesection.\@arabic\c@figure}
\def\fps@figure{h, t}
\@addtoreset{table}{section}
\def\thetable{\thesection.\@arabic\c@table}
\def\fps@table{h, t}
\@addtoreset{equation}{section}

\makeatother

\usepackage{pslatex}
\author{ Jair Koiller\footnote{On a CAPES-Fulbright visit to Caltech, Winter 2005-2006.}\\
Funda\c c\~ao Getulio Vargas \\ Praia de Botafogo 190 \\ Rio de Janeiro, RJ, 22250-040, Brazil  (jkoiller@fgv.br) \\
Kurt Ehlers  \\ Truckee Meadows Community College \\
7000 Dandino Boulevard \\ Reno, NV, 89512-3999, USA (kehlers@tmcc.edu) }

   \title{Rubber rolling over a sphere}  
 
\begin{document} 

\newcommand{\R}{\hbox{\bf R}}

\newtheorem{proposition}{Proposition}
\newtheorem{theorem}{Theorem}
\newtheorem{lemma}{Lemma}
\newtheorem{definition}{Definition}
\newtheorem{notation}{Notation}
\newtheorem{remark}{Remark}
\newtheorem{assertion}{Assertion} 
\newtheorem{conjecture}{Conjecture}
\newtheorem{corollary}{Corollary}
\newtheorem{exercise}{Exercise}
\newtheorem{problem}{Problem}
\newtheorem{criterion}{Criterion}

\addtolength{\topmargin}{-.775in}
\addtolength{\textheight}{1.5in}


\newcommand{\vU}{\vec{U}}
\newcommand{\vX}{\vec{X}}
\newcommand{\vx}{\vec{x}}
\newcommand{\vV}{\vec{V}}
\newcommand{\vv}{\vec{v}}
\newcommand{\vn}{\vec{n}}
\newcommand{\vN}{\vec{N}}
\newcommand{\vu}{\vec{u}}
\newcommand{\vT}{\vec{T}}
\newcommand{\vF}{\vec{F}}
\newcommand{\vW}{\vec{W}}
\newcommand{\vc}{\vec{c}}
\newcommand{\vr}{\vec{ r}}
\newcommand{\vd}{\vec{ d}}
\newcommand{\va}{\vec{ a}}
\newcommand{\cd}{{\cal D}}
\newcommand{\ta}{\tilde{a}}
\newcommand{\tb}{\tilde{b}}
\newcommand{\tc}{\tilde{c}}
\newcommand{\T}{{\bf T}}
\newcommand{\ho}{\hat{\Omega}}
\newcommand{\hT}{\hat{T}}

\maketitle

\pagestyle{headings}

\bigskip
\noindent {\it Key words: Nonholonomic mechanics, Reduction, Chaplygin systems} \\

\noindent {\it AMS MSC(2000) }; 37J60, 70F25, 58A15, 58A30 .

\bigskip

\noindent \small{Abstract.     ``Rubber'' coated rolling bodies   
 satisfy a no-twist in addition to the no slip  satisfied by 
  ``marble'' coated  bodies \cite{EhlersIUTAM}.  Rubber rolling has an interesting differential
  geometric appeal because the geodesic curvatures of the
  curves on the surfaces at the  corresponding points are equal.  The associated distribution in the 5 dimensional configuration space
has 2-3-5 growth (these distributions were first studied by Cartan; he showed that the maximal symmetries occurs  for rubber rolling 
of spheres with  3:1 diameters ratio and materialize the exceptional
group $G_2$).  The 2-3-5  nonholonomic geometries are classified in a companion paper \cite{EhlersKoiller} via 
Cartan's equivalence method \cite{Cartan1}.
  Rubber rolling of a convex body  {\it over a sphere}  defines a generalized Chaplygin system \cite{Chaplygin0,Stanchenko,Iliyev,Kurt2,KoillerArma}
with $SO(3)$
 symmetry group, total space $Q = SO(3) \times S^2$ and base $S^2$,
that can be
reduced to  an almost Hamiltonian system in  $T^*S^2$ with a non-closed 2-form $\omega_{NH}$. In this paper we present some basic results on
the sphere-sphere problem: a dynamically asymmetric but balanced  sphere of radius $b$
(unequal moments of inertia $I_j$ but with center of gravity at the geometric center),  rubber rolling over another sphere of radius $a$.  
  In this example   $\omega_{NH}$  is  conformally symplectic \cite{Vaisman}:
the reduced system becomes Hamiltonian after a
coordinate dependent change of time. In particular there is an invariant measure, whose density  is
the determinant of 
the reduced Legendre transform,  to the power  $ p = \frac{1}{2} ( \frac{b}{a} - 1) $. 
 Using sphero-conical coordinates we verify the results by Borisov and Mamaev \cite{BM1,personalcommunication} that the system is integrable for $p = -1/2$ (ball over a plane)
 and  $p = -3/2 $ (rolling ball with twice the radius of a fixed internal ball).

\newpage
\tableofcontents

\newpage

\section{Introduction } 

   \indent \mbox{    }  Given a  riemannian manifold $Q^n$ 
and a  $s<n$  dimensional distribution $ {\cal D}$  of  subspaces of  $TQ$,
two different theories  apply. One  is  
{\it subriemannian geometry} \cite{Montgomery}    which  appears   in the study of  underactuated  control systems  and gauge theories.  
{\it Nonholonomic geometry} \cite{Arnoldetal,CortesManolo, NF, Bloch,Cushmanbook}  follows
  d'Alembert's principle, and describes mechanical systems with nonintegrable
constraints\footnote{ Hertz \cite{Hertz} introduced the word  {\it nonholonomic} to designate 
nonintegrable distributions and was the first to call attention that, although using the
 same ingredients, the two theories are quite different.}.

This an a companion
 \cite{EhlersKoiller} paper belong to nonholonomic geometry.  We use the following terminology: ``Rubber'' coated bodies means adding the no-twist condition to the usual no-slip constraints of  ``marble'' coated bodies. The dynamics for the former lives on a 7 dimensional phase space, 8 dimensional for the latter.
``Waxed'' marble bodies 
have no velocity constraints: skiding dynamics lives on a 10 dimensional phase space.

Strangelly enough, while marble  bodies have been extensively studied in the NH-literature,
rubber rolling seems to have been neglected. This is curious because they should be easier to study and
furthermore, rubber rolling has an appealing differential geometric interpretation.

 In this paper we provide details of some results on rubber rolling over a sphere announced in \cite{EhlersIUTAM}. Details about the classification of 2-3-5 nonholonomic
geometries via Cartan's equivalence \cite{Cartan1} will be submitted elsewhere.
The expert can go directly to section 3, 
or pass quickly over section 2. Although we tried to be reasonably self-contained, we use freely the geometric mechanics jargon, see eg.
\cite{Arnold,Arnoldetal,MR,Oliva}. We also presented an extensive but focused
list of  references\footnote{We apologize for omissions, specially on nonholonomic
reduction and special geometric structures (like almost Poisson and almost Dirac); we just mention some new connections of nonholonomic systems  in robotics and control \cite{Schaft},\cite{Yoshimura}.}.

 Using sphero-conical coordinates we verify the results of Borisov and Mamaev \cite{BM1,personalcommunication} that the system is integrable for $b/a = 0$ 
 and $b/a=-2$.  A table of results  similar to that of Borisov and Mamaev (tables 1 and  2  of \cite{BM})  for marble bodies rolling over a 
plane or a sphere) is in order. 

\subsection{ Extrinsic differential geometry of rubber rolling}  \label{rollingsection}

In this section we discuss the implications of  rolling  a surface $\Sigma_2$ 
without sliping or twisting  over  a surface $\Sigma_1$  
from  the {\it extrinsic differential geometric} perspective.
  In   \cite{EhlersKoiller}  we show that this process belongs actually
to their {\it intrinsic} geometries.

Given a   curve  $C $  in a surface $\Sigma$,  define the $ (C,\Sigma)$ {\it 
adapted frame}  as
the moving  frame  $F = (t,u, N)$, where $t$ is the tangent vector, $N$ the
 normal to the surface  (so an orientation is chosen) and $u = N \times t $ is the surface
 normal to the curve.   Denoting  $\,\, ' = d/ds$  the derivative with respect to arc length, as it is well known we have
\begin{equation} \label{streq}  t' = \kappa_g \,u + \kappa_n \, N\,,\,\, u' = - \kappa_g \,t - \tau_g \, N\,,\,\,N' = - \kappa_n \,t + \tau_g  \, u\,
\end{equation}
where  $\kappa_g = t' \cdot u  \,$ is called the {\it geodesic curvature}, $\, \kappa_n = t'\cdot N \,\,$
is the {\it normal curvature},  and $\, \tau_g = N' \cdot u  \,\,$  the {\it geodesic torsion}.
Recall an elementary result from classical differential geometry  
\begin{proposition} \label{propstruik} (see e.g. Struik,  \cite{Struik}, p.201, exercise 2, section 4-8.)
The geodesic curvature of a curve $C$ on a surface $S$ is equal
to the ordinary curvature of the plane curve into which $C$ is deformed when the
developable surface enveloped by the tangent planes to $S$ along $C$ is rolled
out on a plane. 
\end{proposition}
\begin{remark}
The inverse problem is hard: finding a curve with prescribed geodesic curvature
on a surface $S$ with metric
$ds^2 = E du^2 + 2F dudv + G dv^2$   gives rise to a nonlinear second order equation for  $u(s),v(s)$. 
\end{remark}
 
Writing  $F$ as the orthonormal matrix
with colums t,u,N,  then $ F^{-1} F' = A$, where  $A$ is the skew symmetric matrix 
$$ A = \left( \begin{array}{lll}   0 & - \kappa_g  & - \kappa_n \\  \kappa_g & 0 &  \tau_g \\ \kappa_n & - \tau_g & 0 
\end{array} \right)  \,\,\,.
$$
Let us write the structure equations (\ref{streq}) as
 \begin{equation}   t' = \omega_F \times t \,,\,\, u' = \omega_F \times u \,,\,\, N' = \omega_F \times N  \,\,\,\,,\,\,\,\,
  \omega_F = - \tau_g \,t - \kappa_n \, u +  \kappa_g \, N  \,\,.
  \end{equation}
This yiels a nice dynamic interpretation: as the curve is traversed with unit velocity,   the  frame is rotating instantaneously
around the vector $\omega_F$ with angular velocity $|| \omega_F || $.

Consider two  surfaces $\Sigma_2$ and  $\Sigma_1$   an suppose  that 
$\Sigma_2$ moves  always  ``touching''  $\Sigma_1$. 
We neglect physical intersections that may occur if one of them is not convex).  At the contact point  the normals will be 
 equal or opposite (depending
on which orientations one choses). 
We use the following conventions: for a closed convex $\Sigma_2$, the exterior normal  defines its orientation;  if $\Sigma_1$ is 
also convex we take the opposite orientations when they roll externally to each other.

When there is no sliding, the tangent vectors at the corresponding curves of points of contact  have the same length, 
so we may use the arc length as a common  parameter for  both curves.
Let 
$C_1:  x = x_1(s)$ and $C_2: x = x_2(s)$ denote the corresponding  curves of points of contact in $\Sigma_1$ and $\Sigma_2$.  
and let us describe the rolling action 
\begin{equation} g(s)= (R(s), a(s))  \in SE(3) \,\,\,
\end{equation}
($g$ acts on $x \in \R^3$ by  $g \cdot x  = Rx + a$, 
meaning
´´first rotate, then translate'').   

Let $\omega_R$ the angular velocity in space of the rotation, and assume the condition 
 $ \omega_R \times x_2 + R^{-1} \dot{a} = 0$.   Then $R \dot{x}_2 = \dot{x}_1 $ 
 and in fact the corresponding adapted   frames are related  by
\begin{equation}   R(s) F_2 (s) = F_1(s)   \,\, ,
\end{equation}
i.e,  they  match under the action of $R(s)$.  Here we  orient  $\Sigma_1$ so that the normals to the surfaces  point in the same direction.  
Differentiating, we get
    $$ F_2^{-1} (R^{-1} R') F_2   +  F_2^{-1} F_2' = F_1^{-1} F_1'
    $$
There is no loss in generality in assuming that at the point of contact corresponding
to $s=0$ we have $F_1(0) = F_2(0) = id$ (so $R(0) = id$, and $x_1(0) = x_2(0) = 0$.   Therefore, if
$ R'(0) $ is the skew-symmetric matrix $B$ given by
$$  B_{12} =  - \omega_3   \,\, , \,\, b_{13} =  \omega_2   \,\,\  ,  \,\, a_{23} =  - \omega_1 \,\,, $$ 
then  
$$   \omega_{F_1} = \omega_{F_2} + \omega_R
$$ 
which gives the following relations\footnote{We learned this result from Mark Levi \cite{Levi}. }  between the invariants of the curves $C_1$ and $C_2$:
\begin{equation}
 \kappa_g^{(1)}  =   \kappa_g^{(2)} +  \omega_3   \,\,\,,\,\,\,\kappa_n^{(1)}  = \kappa_n^{(2)} -  \omega_2  \,\,\,\,,\,\,\,\,
 \tau_g^{(1)}  =  \tau_g^{(2)} -  \omega_1   
\end{equation}


A very important consequence for our purposes  is the following
\begin{proposition} \label{geosame}
 Under the condition of
no twisting ($\omega_3 = 0$) the {\it geodesic curvatures of the contact curves
at the corresponding points are the same.}
\end{proposition}

\subsection{Main results}

 
 Rubber rolling of a {\it convex } body ${\cal B}$ {\it over a sphere of radius $a$} yields a  {\it 
generalized Chaplygin system},  
see \cite{Kurt2,Fedorov,KoillerArma,Stanchenko}.   We have a principal bundle with total space $Q=SO(3) \times S^2$ and base space   $B=S^2$.   
The symmetry group $G$ is $SO(3)$,  acting diagonally. 
The metric is left invariant,
and the constraints define a connection on the principal bundle, called
the ``rubber connection''.   

 The theory developed in \cite{Kurt2} implies that the dynamics reduces to $T^*S^2$, with   a
non-closed 2-form  $\omega_{nh} = \omega_{can}^{T^*S^2} + (J,K)$ and 
 a ``compressed'' Hamiltonian.   The term $(J,K)$ is a semi-basic form, where $J$ is the momentum map and $K$ is the
curvature of the connection.  The configuration variable in the base, $\gamma \in S^2$, has a nice geometric 
interpretation. It is the {\it Poisson vector}, namely, {\it minus} the external normal vector of
the rolling body at the contact point (or equivalently, the normal vector
to the base sphere or plane), seen in the body frame.

If a solution $\gamma(t)$ is found, then the angular velocity $\Omega(t)$ of ${\cal B}$ in the body frame  is the horizontal lift of
$\dot{\gamma}$ via the rubber connection.  To reconstruct the  attitude matrix  one needs to solve a linear system with time dependent coefficients,
$\dot{R} = R [\Omega(t)]$. 
To complete the reconstruction, note that the contact point in $S_a$ is  $q(t) = a R \gamma(t)$. 

In principle, the curve $q(t)$ in $S_a$ could also be reconstructed from $Q(t) =  dN_{\cal B}^{-1}(- \gamma(t))$ on 
$S_b$ (where $dN_{\cal B}$ is the Gauss map of ${\cal B}) $ using the fact that {\it the  geodesics curvatures $\kappa_g$ of $q(t)$
and $Q(t)$  at corresponding
points are the same}.  If one succeeds to solve this nonlinear equation, then the rotation matrix $R$ can be computed, 
by noting that it takes the adapted frame
$(\dot{Q}/|\dot{Q}, \dot{Q} \times \gamma/|\dot{Q} \times \gamma|, -\gamma)$ along $Q(t)$
to the corresponding adapted frame 
$ (\dot{q}/|\dot{q}, q/a \times \dot{q}/|\dot{q} \, ,  \, q/a \,) $,
(modulo a trivial change of signs depending on the chosen orientations).

In passing, we observe that while it is simple to reconstruct a curve in the plane from its
curvature $\kappa(s)$, for a general surface this task results on a complicated nonlinear equation (for $\kappa_g=0$ this is already the geodesic equations). For a spherical
curve we get a linear system with variable coefficients\footnote{Is this solvable by quadratures?}
$ \frac{d^2 q}{ds^2} = \frac{\kappa_g(s)}{a} \,q \times \frac{dq}{ds} - \frac{1}{a^2} q  \,\,.
$


In this paper we discuss the  example where  the surface of the rolling body ${\cal B}$ is  a also a sphere,  of radius $b$. 
It is a ``Chaplygin  sphere'', meaning a sphere  of mass $\mu$  where the center of mass is the geometric center, but dynamically asymmetric, i.e.,  
unequal moments of inertia $I_j$.  The metric is given by 
\begin{equation}
  2T = \mu (1 \pm b/a)^2 || \dot{q} ||^2  +  (A\Omega, \Omega) \,\,\,,\,\,\, A = {\rm diag} I_1, I_2, I_3   
\end{equation}
where $q \in S_a$ is the contact point and $\Omega, \omega$ are respectively the angular velocity of
the rolling sphere $S_b$ with respect to body and space frames, respectively. The  constraints are given by
\begin{equation}
 (1 \pm b/a) \dot{q} = ( \pm b/a) \omega \times q\,\,\,,\,\,\,  \dot{q} \cdot \omega = 0   \,\,.
\end{equation}
(the plus sign corresponds to the external case). 

Hence, when the rolling body is also a  sphere,  the  constraints are invariant under the {\it  right } action of $SO(3)$ 
on the second factor of  $S^2 \times SO(3)$. This system is (morally speaking) akin to
a  {\it  LR nonholonomic Chaplygin  system }  \cite{Fedorov}.    For a true  LR system with 2 dimensional base space one can
 guarantee the existence of a function $f: S^2 \rightarrow \R$ such that  $ d(f \omega_{nh}) = 0 $.
Moreover,  the conformal factor $f$
would be  $f = F^{-1/2}$,  the inverse of
the square root of the determinant of the reduced Legendre transform  (identifying
$TS^2 \equiv T^*S^2$),     where
\begin{equation}
F(\gamma)  = {\rm det}(Leg_{red}) = (I_1I_2I_3)(1 + \frac{b}{a})^2 \, \left( (A^{-1} \gamma, \gamma) + \mu\,b^2 [ \frac{\gamma_2^2 + \gamma_3^2}{I_2I_3} + 
\frac{\gamma_1^2 + \gamma_3^2}{I_1I_3} + \frac{\gamma_1^2 + \gamma_2^2}{I_1I_2}] + \frac{\mu^2 b^4}{I_1I_2I_3} \right) \,\,.
\end{equation} 

Indeed, we will show that in the rubber sphere-sphere problem,   the reduced system is  Hamiltonizable.   However,  we get a different exponent:   
\begin{equation} f_{a,b}(\gamma) = F^{\frac{b-a}{2a}}\,,\,\ \gamma \in S^2 \, .
\end{equation}

Probably an explanation for this mysterious exponent will come from a  
study of  Chaplygin systems  of the form $G  \hookrightarrow G \times H/G  \rightarrow H/G$
where the base space is a homogeneous space, and $G$ acts diagonally
in the total space. 
This theory should build up fom the   LR systems  $H  \hookrightarrow G    \rightarrow H/G$ studied by Fedorov and Jovanovic \cite{Fedorov}.


Borisov and Mamaev  \cite{BM1} observed that  by taking suitable combinations
 of the reduced variables,  solutions of 
 rubber rolling of a sphere over a plane can be mapped to the solutions of
  Veselova's system \cite{Veselovs,Veselovs1}.   They have also shown  integrability \cite{personalcommunication}
(in this volume)
in the case $b=-2a$,  where the fixed ball has half the radius and is internal to the rolling ball.  We confirm their result by showing that the
 sphero-conical coordinates separates the Hamiltonian in the new time.  
 Using sphero-conical coordinates$I_1 < \lambda_1 < I_2  < \lambda_2 < I_3$, given by
\begin{equation} 
\left(
 \gamma_1^2 \,,\, \gamma_2^2 \,,\,  \gamma_3^2 
\right ) =
\left (   \frac{(I_1 - \lambda_1)(I_1 - \lambda_2)}{(I_1 - I_3)(I_1 - I_2)}  \,,\, \frac{(I_2 - \lambda_1)(I_2 - \lambda_2)}{(I_2 - I_3)(I_2 - I_1)}  \,,\, \frac{(I_3 - \lambda_1)(I_3 - \lambda_2)}{(I_3 - I_1)(I_3 - I_2)} 
  \right )  \,\,\, .
\end{equation}

Quantitatively, our main result can be given as follows: 
\begin{theorem} Let $(q,R) \in S^2_a \times SO(3)$ the configuration space
coordinates, where $q$ is the contact point in the base sphere $S_a$ and
$R$ the attitude matrix of the moving sphere $S_b$.  Let $ \gamma = R^{-1} (q/a) =
- Q/b$ be
the Poisson vector, where $Q$ is the contact point in the moving sphere,
seen in the body frame.
 The equations of motion for $\gamma$ are governed, in a new time scale $\tau$ 
such that
\begin{equation}  d\tau/dt = F(\gamma)^{\frac{b-a}{2a}} \,\,\,,\,\,\, F(\gamma) = (1+\frac{b}{a})^2 \, (\lambda_1 + \mu b^2) (\lambda_2 + \mu b^2)  \,\,\,
\end{equation}
by a Hamiltonian system 
 (the  nonholonomic vectorfield is  $X_{nh} = F(\gamma)^{\frac{b-a}{2a}} X_H$)
\begin{equation}  \label{hamiltoniansph}
X_H:\,\, \frac{d\lambda}{d\tau} = H_P\,\,,\,\,  \frac{dP}{d\tau} = - H_{\lambda} \,\,\,,\,\,\,\,  
2H =  \frac{P_1^2}{[(\lambda_1+\mu b^2)(\lambda_2+\mu b^2)]^{\frac{b-a}{a}} c_1} +  \frac{P_2^2}{[(\lambda_1+\mu b^2)(\lambda_2+\mu b^2)]^{\frac{b-a}{a}}c_2} 
\end{equation}
where
\begin{equation}
c_1 = \frac{1}{4}(1+ \frac{b}{a})^2 (\lambda_2-\lambda_1)  \frac{\lambda_2 + \mu b^2}{(\lambda_1-I_1)(I_2-\lambda_1)(I_3-\lambda_1)}\,\,,\,\, 
c_2 = \frac{1}{4}(1+ \frac{b}{a})^2 (\lambda_2-\lambda_1)  \frac{\lambda_1 + \mu b^2}{(\lambda_2-I_1)(\lambda_2-I_2)(I_3-\lambda_2)} \,\,\,.
\end{equation}

\end{theorem}
\begin{corollary}
The terms $\lambda_2+\mu b^2$  
in $c_1$ and  $\lambda_1+\mu b^2$ 
in $c_2$ are a nuisance for
Hamilton-Jacobi separation, but  they disappear in two cases, discovered first
by Borisov and Mamaev \cite{personalcommunication}.  One is seen immediately, $(b-a)/a = -1 $, that is,
$b/a \rightarrow 0$, the planar case. The other is   $(b-a)/a = -3 $, where we have the cross factors $(\lambda_2 + mu b^2)/(\lambda_1 - \lambda_2)$
in  $P_1^2$ and $(\lambda_1 + mu b^2)/(\lambda_1 - \lambda_2)$ in  $P_2^2$; 
to see that that the Hamiltonian also separates in this case, multiplying both sides of (\ref{hamiltoniansph}) by
$$   \frac{\lambda_2 - \lambda_1}{(\lambda_1 + mu b^2)(\lambda_2 + mu b^2)} =   \frac{1}{\lambda_1 + mu b^2} -
  \frac{1}{\lambda_2 + mu b^2}
$$

It is not known if there are hidden
integrals for other values of $b/a$ and inertias $I_j$ (except for the case of 
two equal inertias), but this is unlikely.  The case $b=a$ is special because the constraints are {\it holonomic}, so  we have directly a  
two degrees of freedom Hamiltonian system.  
It was somewhat frustrating
to realize that in this case the problem  {\it does not} separate
in sphero-conical coordinates  and in fact seems to be chaotic from numerical experiments\footnote{We admit losing a bet to Ivan Mamaev
and Alexei Borisov, and we owe them a dinner in a barbecue house in Rio.}.
\end{corollary}

\section{Preliminaries}

\subsection{Skiding dynamics} 

\indent \mbox{} As we all know, the boundary of a {\it strictly convex}  body  ${\cal B}$  is a closed surface $ \Sigma_2 $ with strictly positive Gaussian curvature at all points. 

We say ${\cal B}$  is in {\it standard position}  when  $0$ is the center of mass and  the principal axis of inertia are aligned with
 OX,OY,OZ with moments of inertia  $I_1 \leq I_2 \leq  I_3$,  respectively.
Under a Euclidian motion $g = (x,R) \in SE(3)$ the surface goes to 
$ g \cdot \Sigma_2 $,  the center of mass of  $ g \cdot {\cal B}$ 
is at  $x \in \R^3$ and
the three principal axis become, respectively,  the columns  $e_1, e_2, e_3$ of the {\it attitude matrix}  $R$.  

The {\it  Gauss mappings}   $N_j : \Sigma_j \rightarrow S^2$ are important ingredients for the sequel. We make  no special requirement on $N_1$, but
 assuming  $\Sigma_2$ to be convex  is useful, it  guarantees  that $N_2$ is one to one and onto $S^2$.

  The reader can easily sketch a figure to visualize our notation.  
 Let $ \Sigma_1$  be  fixed and $\Sigma_2$  move under the action of  $SE(3)$.
The unconstrained configuration space is {\it ten} dimensional, 
$ (SE(3) \times \Sigma_2) \times \Sigma_1 \,\,.$  
When we impose the condition that  the two bodies {\it  touch},  the configuration space becomes {\it  five} dimensional. 
One way to do the dimension count is to consider  equations
$F_1(x_1,x_2,x_3) = 0 $ for  $\Sigma_1$ and $F_2(Q_1,Q_2,Q_3) = 0$ for  $\Sigma_2$.  The {\it touching manifold}  
$M^5 \subset  (SE(3) \times \Sigma_2) \times \Sigma_1$
consists of solutions for the 7 equations in 12 variables
$$ F_1(X) = 0 \,\,,\,\,  F_2(Q) = 0 \,\,, g. Q = X\,\,, \,\,{\rm and}\,\,\,\,  \nabla F_1 (X)/|F_1 (X)| = \nabla  F_2( g^{-1} Q ) \,\, . $$
We have 5 effective  equations as the last three count as two.  
In the sequel we will identify  $M^5  \equiv  \Sigma_1 \times SO(3) \,, $
 where the first component is the contact point, and the second component the attitude matrix of the rolling body.    
Given $q_1 \in \Sigma_1$ and $R \in SO(3)$ we can find 
the corresponding point $Q_2 \in \Sigma_2$  and  the position
$x = x_{CM}$ of the center of mass.  Indeed, since the normals are aligned, we have
$ N_1(q_1) = - R N_2(Q_2)  \,\,\,, $
where we orient $\Sigma_2$ with the exterior normal. For visualization convenience, we orient   $\Sigma_1$ in such a way that 
at the contact points the normals point in opposite ways (for instance, in  the sphere-sphere case, both are the exterior normals).  Hence
\begin{equation}  \label{Q2}
   Q_2 = N_2^{-1} \, (- R^{-1}\, N_1(q_1) )
\end{equation}
\begin{equation}    x = x_{CM} = q_1 - R Q_2 = q_1 - R \, ( N_2^{-1} ( - R^{-1} N_1(q_1) ) ) \label{CM}
\end{equation}
(in the sequel we may drop the suffixes when no confusion may arise).


As found long ago by Euler, the kinetic energy of a rigid body  ${\cal B}$ is given by
\begin{equation}  2 T(x, R, \dot{x}, \dot{R})  =  \mu || \dot{x} ||^2 +   I_1 \Omega_1^2 + I_2 \Omega_2^2 + I_3 \Omega_3^2  \,\,\,,  \label{kinetic}
\end{equation}
where
 $\mu$ is the total mass, $I_j$ the moments of  inertia about the principal axis    $e_1, e_2, e_3$  attached at the center of mass
  $x$, and  $ \Omega = (\Omega_1, \Omega_2 , \Omega_3 )^t $  the  angular velocity vector written  {\it the body frame}.  
More precisely,  $R^{-1} \dot{R} = [\Omega] $ is the skew symmetric matrix
$$ [ \Omega ] = \left( \begin{array}{lll}  0  & - \Omega_3  & \Omega_2 \\
\Omega_3 & 0 & - \Omega_1 \\ 
- \Omega_2 & \Omega_1 & 0  
\end{array}   \right)
$$
See Arnold \cite{Arnold} for lower/upper  case notations:
 {\it The angular velocity in space} is    $\dot{R} R^{-1} = [\omega ]$,  with $ \omega = R \Omega $.

There is no coupling between the translational and the rotational motions if the body is moving freely in space. 

What is the dynamics when we impose the (holonomic) constraint of skiding?
For simplicity, we ignore potential forces. In order to  compute  the Lagrangian
 \begin{equation}  \label{lagrangianofQ}
 T(q_1,R,\dot{q}_1,\dot{R}) \,\,\, \, {\rm  on}  \,\,\,  Q = \Sigma_1 \times SO(3) 
\end{equation}
 governing  the {\it holonomical system} defined  by   the touching conditions  (\ref{Q2},\ref{CM}),   we replace  in (\ref{kinetic}) the velocity $\dot{x}$ of the center of mass by
\begin{equation} \dot{x} = \dot{q}_1 - D_{(R,q_1)} \, [R \, ( N_2^{-1} ( - R^{-1} N_1(q_1) ) )] \,(\dot{R},\dot{q}_1 ) \,\,.  \label{dotx}
\end{equation}
This task is  not so simple, because we need  the derivative of the map
$$ F: (R,q_1) \in SO(3) \times \Sigma_1 \mapsto  R \, ( N_2^{-1} ( - R^{-1} N_1(q_1) ) ) \in \R^3 $$
appearing in (\ref{dotx}). This notwithstanding,  all the geometric information  we need  is contained in the Gauss maps:
\begin{lemma} \label{gaussmaps}
\begin{equation}
 dF_{|(R,q_1)} \, \cdot (\dot{R}, \dot{q}_1) = \omega \times (q_1 - x) + (dN_{g \Sigma_2})^{-1} \,( \omega \times N_1 (q_1) - dN_1 (q_1) \cdot \dot{q}_1 )
\end{equation}
\end{lemma}
\noindent {\bf Proof.} This is a simple exercise on Advanced Calculus.
Here the
rotational velocity $\dot{R}$ is written  in  the space frame,
$ \dot{R} R^{-1} = [\omega]$. Denoting by $ g(t) \cdot \Sigma_2$   the current (located) position of the moving body, 
 $ g = (R,x) \in SE(3)$,  we can  reinterpret some objects that appear
in the derivation,  for instance
$    R^{-1} dN_2 (Q_2) R = d N_{g \Sigma_2}(q_1)     \,\,.
$

\noindent {\it Example: sphere-sphere skiding.}  Let  $\Sigma_1 = S_a$
a sphere of radius $a$ centered at the origin, $ \Sigma_2 = S_b$  a sphere of radius $b$ rolling over $\Sigma_1$. 
 Then (\ref{CM}) becomes as expected to $ x = (1 \pm b/a) \, q_1 $
  The plus sign corresponds to the external case. The Lagrangian is
\begin{equation}  \label{spherelagrangian}
  2T = \mu (1 \pm b/a)^2 || \dot{q}_1 ||^2  +   I_1 \Omega_1^2 + I_2 \Omega_2^2 + I_3 \Omega_3^2
\end{equation}
and the dynamics uncouple.  Solutions are  great circles in the sphere and Euler rigid body motion for the moving ball.

\begin{proposition}  \label{skidingonsphere}
For the problem of a convex body ${\cal B}$ with surface $\Sigma_2$ skiding over  a sphere  $\Sigma_1 = S_a$  there is $SO(3)$  symmetry.
 We have a principal bundle structure  
with total space  $Q = S_a \times SO(3)$,  
\begin{equation}  \label{bundle1}
 SO(3) \hookrightarrow   Q  \rightarrow  Q/SO(3) \,\,\,,\,\,\,  R \cdot  (q, S) = (R(q), RS ) 
\end{equation}
given  by the diagonal action. 
\end{proposition}
Along the fiber  ${\cal B}$ is been moved rigidly    rotations  around the {\it center} of sphere $S_a$.  In a sense motion along the fiber is ``pure'' skiding.
The base of the bundle is the sphere $S_a$  and we can
normalize to the unit sphere, parametrized by $\gamma$.
There is  a global section,
namely  
$  q_1 = a \gamma \in S^2 \mapsto (a \gamma, I) \in S_a \times SO(3) \,\,  .$
Notice that the induced  global trivialization  is {\it  not} given just by the direct product structure.   
Fibers are
diffeomorphic to $SO(3)$ via $ R \mapsto (a R(\gamma), R)$. Along the global
section,   $q_1 = a \gamma $
is the contact point in $\Sigma_1 = S_a$, $ Q_2 = N_2^{-1} (-\gamma)$ is the corresponding contact point in $\Sigma_2$ (in its standard position). 
Given  a curve $R(t) \in SO(3), \,\, R(0) = I,\,  \dot{R}(0) = [\rho]$,
\begin{equation}  V(\rho)  = \frac{d}{dt}_{|t=0} \, R(t) ( a \gamma, I) = ( a \rho \times \gamma, \rho) 
\end{equation}
gives a vertical vector. 

\begin{remark}
Reduced equations for the skiding dynamics can be obtained by Marsden-Weinstein procedure \cite{MW}.  For the symplectic reduction we need 
  the momentum mapping  $J:T^*Q \rightarrow \R^3$.   An element of $T^*Q$ is of the form $(p_R,p_{q_1})$, where 
 $p_R \mapsto  m \in sO(3)^* \equiv \R^3 $ via right translation to the identity,  and $   p_{q_1} \in  T^*(S_a)$ can be represented 
by a vector $r$  perpendicular  to $q_1$.
Using the   abstract nonsense  rule   
 `` $\,\, J \,\,{\rm of}\,\,  p \,\, {\rm at}\,\,\, q  \,\,, \,\,\,  {\rm on} \,\,  X \,\, \, = \,\, \, p\,\, {\rm at} \, \, q \,\, {\rm on } \,\,  V(X) \, $ `` 
we get  
\begin{equation}  \label{momentummap}
J: T^*(S_a  \times SO(3)) \rightarrow  sO(3)^* \, (\equiv \R^3)\,\,\,\,,\,\,\,\,  J(p_R, p_q) = m + q_1 \times r = {\ell}  \,\, .
\end{equation}
The  six dimensional reduced
symplectic manifolds  are $M^6 = J^{-1}(\ell)/S^1$,  where
 $S^1$ denotes the isotropy group of $\ell$.   
 \end{remark}

\subsection{ No slip and no-twist constraints:  kinematical relations for the Poisson vector} 

The  {\it  Poisson vector}
\begin{equation}  \label{Poissonvector}
\gamma = R^{-1} N_1(q_1) = - N_2(Q_2) 
\end{equation}
is the normal vector to $\Sigma_1$ seen in the body frame and is the basic object for rubber rolling.  It is related to $Q_2$ by {\it minus} the Gauss map. 
We  get immediatelly from  $R \dot{Q}_2 = \dot{q}_1 $:
\begin{lemma}  Kinematical relations for the Poisson vector, no slip:
\begin{equation} \label{kinematicalrelation}
\dot{\gamma} +  \Omega \times \gamma =   R^{-1} dN_1(q_1) \cdot \dot{q}_1 =  [ R^{-1} dN_1(q_1) R ] \cdot \dot{Q}_2  \,\,\,\,  .
\end{equation}
\end{lemma}
In this formula,$\,\,\, q_1 =  (N_1)^{-1} (R \gamma) $ when the local inversion 
is possible\footnote{Points on $\Sigma_1$ where one of the principal curvatures vanish may be specially  relevant for the dynamics.}. 
In particular, if  $\Sigma_1$ is a sphere of radius $a$ or a plane ($a = \infty$),  the term  $ R^{-1} dN_1(q_1) R$ is just  $1/a$ times  the
identity, 
 the kinematical relation becomes
\begin{equation} \label{bodysph}  \dot{\gamma} +  \Omega \times \gamma = R^{-1} (\dot{q}_1/a) = \frac{1}{a}  \dot{Q}_2  = -  \frac{1}{a} \, d (N_2)^{-1}_{\gamma} \cdot  \dot{\gamma} 
\end{equation}
Furthermore, if $\Sigma_2$ is also a sphere, of radius $b > 0$,
\begin{equation}  \label{kinsphsph}
\dot{\gamma} + \kappa \Omega \times \gamma = 0\,\,\,\,,\,\,\, \kappa = \frac{a}{a \pm b}
\end{equation}
where the plus sign corresponds to the external case, 
minus when one of the spheres contain the other. We can use only the plus sign and allow $b$ to be 
negative to represent the internal case (if $b<0$ and $|b|>a$ then the fixed sphere is inside the rolling sphere).

\subsection{Rubber connection and curvature}

Imposing the no-slip constraints produces a  map 
\begin{equation} \label{rolling}  \omega \mapsto  \dot{q}_1   \,\,\, ({\rm at} \,\,\ q_1)
\end{equation}
which  gives a distribution of  3-subspaces in $Q$.  
To rule out twisting, one adds  
\begin{equation}   \label{notwist}
 \omega \cdot N_1(q_1) = 0   \,\,\, {\rm or} \,\, {\rm equivalently} \,\,\,\, \Omega \cdot \gamma = 0 \,.
\end{equation}   
In more detail, the no-slip constraint  in $Q = \Sigma_1 \times SO(3)$  follows by replacing the left-hand side $\dot{x}$ in
(\ref{dotx})   by  
\begin{equation} \label{dotxx}  \dot{x} = \omega \times (x - q_1)  \,\,. 
\end{equation} 
In view of lemma \ref{gaussmaps}, the result is 
\begin{equation}  \label{rolling1}
\dot{q}_1 = dN_{g \Sigma_2}^{-1} \,(  \omega \times N_1 (q_1) - dN_1 (q_1) \cdot \dot{q}_1 )
\end{equation}
or,  in a more symmetric form,
\begin{equation}  \label{rolling2}
[ dN_{g \Sigma_2} + dN_1 ](q_1) \cdot \dot{q}_1   =    \omega \times N_1 (q_1)   =  \omega \times (- R N_2(Q_2))  \,\, .
\end{equation}
Since for no-slip the tangent vectors of curves  $q_1(t)$ and $Q_2(t)$ correspond,
 $   R \dot{Q}_2 = \dot{q}_1 $, equation (\ref{rolling2}) gives
$$  \Omega \times \gamma =  R^{-1} [ dN_{g \Sigma_2} + dN_1 ](q_1) \cdot \dot{q}_1 = dN_2(Q_2) \cdot \dot{Q}_2 +  dN_1 (q_1) \cdot \dot{q}_1\,\,.
$$

Imposing in addition to (\ref{rolling2}) the no-twist constraint (\ref{notwist})
allows us to solve for $\omega$ in terms of $\dot{q}_1$.
\begin{lemma}  The ``rubber connection'':
\begin{equation}  \label{connec1}
\omega = N_1(q_1) \times ([ dN_{g \Sigma_2} + dN_1 ](q_1) \cdot \dot{q}_1  \,\,\, .
\end{equation}
\end{lemma}
This formula defines an {\it Ehresmann connection} on the bundle
$ SO(3) \times \Sigma_1  \rightarrow \Sigma_1$ . Equation (\ref{connec1}) can
be rewritten in the body frame as
\begin{equation} \label{connect20}
\Omega = [ ( Id + (R^{-1} dN_1 R)(dN_2)^{-1} )] (\dot{\gamma})  \times \gamma
\end{equation} 
 Note that there is no $SO(3)$ equivariance
unless $\Sigma_1$ is a sphere. 

Consider a small curve $C_1$ around a point $q_1 \in \Sigma_1$. The curvature of the Ehresmann connection is the limit
\begin{center}
$\lim \,\,\,$ ( how much a frame attached to $\Sigma_2$
rotates )  $\,\,\, /   \,\,\,\,$  ( area inside the curve in $\Sigma_1$  )  $\,\,\,\,\,\,\,\,\,\,, \,\,\,\,\,$  as  $C_1$  shrinks . 
\end{center}

We have computed the curvature  for the case of a sphere rolling over another sphere, see  section \ref{curvsphsph} below. In order to appreciate the
difficulties of the calculation, we challenge the reader to attempt computing by brute force the curvature
in the case of a surface of revolution rolling over a sphere.
We claim, however, that it is  possible to give a formula for the Ehresmann curvature  of the rubber connection in the case of {\it general surfaces}  $\Sigma_1,\, \Sigma_2$ 
in terms of natural geometric objects.  Accepting at face value the information that rubber rolling belongs to the intrinsic geometry \cite{Bryant}, the idea is to use  proposition \ref{geosame} together   with the
Gauss-Bonnet theorem. The result (see details in the companion paper \cite{EhlersKoiller}) is a kind of dynamical ``Egregium theorem'':
\begin{theorem} (Bryant and Hsu \cite{Bryant}).  Curvature of Ehresmann connection, general case:
Let $k_i(q_i)$ the Gauss curvature of $\Sigma_i$ at the corresponding 
points $q_i$, such that $\, q_1 = g q_2\,\,,\,\,\, g = (R,a) \in SE(3)$. Then
\begin{equation}  \label{Ehrcurv}
K = ( 1 - \frac{k_1(q_1)}{k_2(q_2)} ) \, k_1(q_1)\, d\Sigma_1   \,\,.
\end{equation} 
where $d\Sigma_1$ is the area form of $\, \Sigma_1$.
\end{theorem}
This result is remarkable because the connection 1-form already involves 
the derivative of the Gauss maps, so at first sight the Ehresmann curvature should involve second derivatives.

\begin{remark}
Koon and Marsden  have shown that a key
ingredient for the equations of motion for nonholonomic systems  is the {\it curvature} of the {\it local} Ehresmann connection 
associated to a splitting of a coordinate system (see \cite{BKMM,KM1}). Both for marble
or rubber rolling this splitting is {\it global}.
See Theorem \ref{reductionchaplygin} below for the case when $\Sigma_1$ is a sphere.
\end{remark}

\subsection{Equations of motion: general physics approach }

In this section we follow Borisov and Mamaev  \cite{BM} closely\footnote{For more abstract approaches, see eg., \cite{KM}, \cite{Marle1},\cite{Marle2}.}.
In addition to the kinematical relation (\ref{kinematicalrelation}), three differential equations for the rolling body can be derived from 
$$ \dot{\ell}_{CM} = \tau_{CM}  \,\,\,, $$
where  $\ell_{CM}\,,\, \tau_{CM}\,$ are respectively the angular momentum
and the torque of the external forces with respect to the center of mass.

 In basic mechanics  textbooks it is
sometimes overlooked that  the above  formula  {\it does not hold in general }   if one takes  $q_1$ as the base point instead of the center of mass $CM$.
Recall that the angular momentum with respect to the contact point $q_1 $   and with respect to the center of mass $x$  are related  by
\begin{equation} \label{angmomentumcont}
  \ell_{q_1} = (x - q_1) \times \mu \dot{x} + \ell_{CM}  \,\,.
\end{equation}
Differentiating, we get
$$  \dot{ \ell_{q_1}} = [ (x - q_1) \times f_{ext} + \tau_{CM} ] + \frac{d}{dt}[(x - q_1)] \times \mu \dot{x}\,\,,
$$ 
and we recognize the term  $(x - q_1) \times f_{ext} + \tau_{CM}$
as the torque about the contact point. The correct  form of the torque equation is
\begin{equation}
 \dot{ \ell_{q_1}} =  \tau_{q_1} +   \frac{d}{dt}[(x - q_1)] \times \mu \dot{x}  \,\,.
\end{equation} 
containing an ammended term, which can be called the {\it dynamical torque},
\begin{equation}   \label{dyntorque}
   \tau_{dyn} =  \frac{d}{dt}[(x - q_1)] \times \mu \dot{x}  \,\,\,.
\end{equation} 
The non-sliping condition implies  $\dot{x} = \omega \times (x - q_1) $, so in the space frame
\begin{equation}  
 \dot{\ell_{q_1}} = \tau_{q_1} +   \frac{d}{dt}[(x - q_1)] \times \mu \omega \times (x - q_1)   \,\,.
\end{equation}
If we go to the body frame, we get\footnote{ This is precisely eqs. (1.1) of \cite{BM} in  our  notation.}
\begin{equation} \label{momentumequation}  \dot{L} + \Omega \times L + (\Omega \times Q_2)  \times (\mu \dot{Q}_2 ) =  T
\end{equation}
where  we recall for clarity, that
$$\ell_{q_1} = R L, \,\,  x - q_1 = R Q_2, \,\,  \omega = R \Omega, \,\, \tau = R T  \,\,.$$

In order to obtain a linear relation between $L$ and $\Omega$  (depending on the Poisson vector $\gamma$), we  rewrite 
the  total energy of the system  as
$$ 2H = 2T =  (\ell_{q_1}, \omega) =  (L, \Omega) \,\,,\,\,\,  {\rm whith}\,\,\,  R L = \ell_{q_1},  $$
where $\ell_{q_1}$ is the angular momentum with respect to the contact point viewed in space,  $L$ its coordinates in the body frame.  
Substituting the no-slipping  constraint in  (\ref{kinetic}), one gets 
\begin{equation}   \label{LofOmega}
L = A \Omega + \mu Q \times (\Omega \times Q)\,\,\,,\,\,\,  A = diag (I_1,I_2,I_3) \,\,.
\end{equation}
$ L =  \frac{ \partial T}{\partial \Omega}  \,\,, $
is the angular momentum of the rolling body with respect to the contact point,  seen in the body frame.
Relation (\ref{LofOmega})
  can be inverted,
\begin{equation} \label{OmegaofL}
   \Omega  =  (A + \mu ||Q||^2 id)^{-1} \,L + \alpha(L,Q) (A + \mu ||Q||^2 id)^{-1}\,(Q)     \,\,\, 
   \end{equation}
where
\begin{equation}
\alpha(L,Q) = \mu  \frac{ (Q, (A + \mu ||Q||^2 id)^{-1} L)}{1 -
       \mu  (Q, (A + \mu ||Q||^2 id)^{-1} Q) }  \,\,.
\end{equation} 
It will be  convenient  to use the shorthand  $\tilde{A} = A + \mu ||Q||^2 id $.

 \begin{proposition}
(Counting dimensions and equations) $\,\,\,\,$  In the no-slip case
 we have an {\it  eight  dimensional phase-space}.  Using  $ (R,\gamma, L) \in SO(3) \times S^2 \times \R^3$ as coordinates (we are assuming that
$q_1 = (N_1)^{-1} (R \gamma)$ can be locally inverted), 
the  dynamics is given by the  momentum equation (\ref{momentumequation}),  the  kinematic relations  (\ref{kinematicalrelation}) and the  attitude
equation  $\dot{R} = R [\Omega]$. 
We will have a closed system of equations since  $\Omega $ is a function of $(L,\gamma)$ via formula  (\ref{OmegaofL}). 
If we add the no-twist constraint, the phase space is {\it seven dimensional}:  one adds the constraint $(\Omega,\gamma) = 0 $, and uses it to
eliminate the unknown torque around the normal.  
\end{proposition}

Reconstruction
of $q_1$  motion in $\Sigma_1$ can be done either via  $\dot{q_1} = R \dot{Q}_2 $, (where as we recall,
$Q$ is related to  $\gamma$ via the Gauss map) or simply inverting $N_1(q_1) = R \gamma$.
 In general, the kinematical relation couples $R$ with the other variables, but  when $\Sigma_1$ is a plane or a sphere,  
 the attitude equation {\it  decouples}  from the equation for $\gamma$ and $ L$,  a consequence of $SO(3)$ symmetry.

\begin{remark} If $x$ and $q_1$ move in such a way that a combination of them with fixed coefficients remains constant, then $\dot{x}$ and $\dot{q}_1$ 
are proportional. Consequently, 
 the ammended  term  in \ref{dyntorque} vanishes identically. For instance, this happens  in the case of
spheres rolling over spheres (or a plane)\footnote{Historical note:  Chaplygin seems to have been the first one to 
discuss   geometrical conditions under which the second term vanishes.   See \cite{NF,Arnoldetal,Kozlov}
for   discussions on this often neglected issue.}. 
\end{remark}

\smallskip

\noindent {\bf Sphere-sphere rubber rolling.}  As the dynamical torque is identically zero\footnote{For  abstract explanations, see  discussions on
nonholonomic Noether theorem in \cite{Arnoldetal,Sniatycki98}.} ,  we are left with
\begin{equation}
\dot{L} + \Omega \times L = \tau \gamma\,\,,\,\, \dot{\gamma} = \kappa \gamma \times \Omega\,\,\,\,,  \,\,\,\,  L = \tilde{A} \Omega + \mu b^2 (\Omega, \gamma) \gamma = \tilde{A} \Omega \,\,,\,\,\,\,\,{\rm with} \,\,\,\,\,(\Omega, \gamma) = 0 \,\,\,
\end{equation}
with $ \kappa = a/(a+b) $. 
The right hand  side $T$ is the external torque about the contact point  in the body frame (zero for the no-slip case, but twisting allowed.) 
  In view of the no-twist constraint $ (\Omega,\gamma) = 0 $ it is immediate that  $\Omega = \tilde{A}^{-1} \, L $, and we get the following expression for
the multiplier:
\begin{equation}  \label{tauequation}
\tau \, (\gamma, \tilde{A}^{-1} \gamma) = (\tilde{A}^{-1} \, L \,  \times L, \,\tilde{A}^{-1} \gamma) + \kappa ( \tilde{A}^{-1} L,  \tilde{A}^{-1} L  \times \gamma )  \,\,.
\end{equation}

The problem of a rubber ball rolling on a plane corresponds  to $\kappa = 1$ and was discussed briefly in \cite{Kurt2}.  We showed
the existence of an invariant
measure
\begin{equation}
\nu = F(\gamma)^{-1/2} dL_1 dL_2 dL_3  d\gamma_1 d\gamma_2 d\gamma_3 
\end{equation}
with
 \begin{equation}
F=  I_1I_2I_3\, \left( (A^{-1} \gamma, \gamma) + \mu\,b^2 [ \frac{\gamma_2^2 + \gamma_3^2}{I_2I_3} + \frac{\gamma_1^2 + \gamma_3^2}{I_1I_3} + 
\frac{\gamma_1^2 + \gamma_2^2}{I_1I_2}] + \frac{\mu^2 b^4}{I_1I_2I_3} \right) \,\,.
\end{equation} 
Borisov and Mamaev \cite{BM1} called the attention that
the solutions of the rubber ball on the plane can be set into correspondence with solutions of Veselova's system \cite{Veselovs,Veselovs1}
namely, a rigid body with one right invariant constraint  (formally, 
 $b=0$).    They used a  trick that seems to go back to Chaplygin, which is really a 
``Columbus egg'':  one  takes suitable linear combinations of $\gamma$ and $L$ and
show that the equations of motion correspond.  Unfortunately,  this clever maneuver only seems to work when $\kappa = 1$.


\subsection{Marble Chaplygin sphere over a sphere or a plane} 

This   section is intended  as 
background material, just for comparison with the rubber rolling problem discussed next.
For more details on what is known about marble rolling of convex bodies over 
a plane or a sphere, see    \cite{BM}. 


 In 1903 Chaplygin \cite{Chaplygin1} integrated using hyperelliptic functions  the following system for   $(L,\gamma) \in \Re^3 \times \Re^3$: 
\begin{equation} \label{L}  \dot{L} +  \Omega \times L = 0 \, \,,
\end{equation}
\begin{equation}  \label{gama}
 \dot{\gamma} +  \Omega \times \gamma = 0 \,\,\,,
\end{equation}
where $\Omega$ is defined by (\ref{OmegaofL}) with $Q =  r \gamma$, namely,
\begin{equation} \label{Omega}
   \Omega(L,\gamma)  =  (A + \mu r^2 id)^{-1} \,L + \alpha(L,\gamma) (A + \mu r^2 id)^{-1}\,(\gamma)     \,\,\, .
\end{equation}
with
\begin{equation} \label{alpha}
\alpha(L,\gamma) = \mu r^2 \frac{ (\gamma, (A + \mu r^2 id)^{-1} L)}{1 -
       \mu r^2 (\gamma, (A + \mu r^2 id)^{-1} \gamma) }  \,\,.
\end{equation} 
They describe the motion of a marble ball of radius $r$ and mass $\mu$, moments of inertia
$I_j$ about the geometric center (which is also the center of mass),  rolling without sliping 
over a plane. Twisting motions (rotations about the vertical) are allowed. $L$ is the angular momentum of the ball with respect to the contact point,
$\Omega$ the angular velocity  and
$\gamma$ the vertical vector, all viewed on a  reference frame attached to the ball.

  This system has four independent integrals 
\begin{equation} f_1 = (L,\Omega)\,,\,\, f_2 =  (L,L)\, , \,\, f_3 = (\gamma, \gamma)\,,  \,\,
f_4 = (L,\gamma) \, 
\end{equation}
where $2 f_1$ is the energy
 and  $f_4$, called the {\it area integral}, is the third spacial component of the angular momentum.
 Topologically, the common level sets of these integrals  are tori;  Chaplygin  showed that there is a smooth invariant measure  $ [ F(\gamma)]^{-1/2} \, d\gamma \, dL    \, $, with
\begin{equation}  \label{chapfunction}
F(\gamma) = \frac{1}{\mu r^2}  -  (\gamma,  \tilde{A}^{-1} \gamma)\,,\,\, \tilde{A} = A + \mu r^2 id 
\end{equation}
Integrability follows from the celebrated {\it Jacobi's last multiplier method}, which says that
  an ODE on a {\it two} dimensional manifold having a smooth invariant measure can be
solved by quadratures. Detailed remakes of Chaplygin's paper appeared recently, see Kilin \cite{Kilin} and   Duistermaat \cite{Duis}. 
 The planar motion can  be reconstructed by integrating
$  \dot{x} = r \omega(t)  \times k $.
\begin{remark}
The issue whether Chaplygin marble system is hamiltonizable or not  has been discussed in \cite{Borisovsphere, Borisovsphere1, Duis}. Perhaps due to some misprints we could not verify the Borisov-Mamaev bracket, but 
we are glad to know that it has been verified  and
explained independently by  Naranjo\cite{Naranjo}.  Theoretical issues about
almost Poisson brackets are outside the scope of our paper, but we believe
will be much in evidence in the next years, see \cite{Cantrijn1,Cantrijn3,Marle2,Marle3} for background.
\end{remark}

To extend  Chaplygin's sphere-plane equations  to the sphere-sphere case, it suficies to replace  equation (\ref{gama})  by  the kinematical relation
\begin{equation} \label{kapa}  \dot{\gamma} + \kappa  \Omega \times \gamma = 0 
\end{equation}
with $\kappa = a/(a + b)$ (Recall that $b>0$ in the external case, so   $0 < \kappa <1$;   in the internal case where $b<0$ and $|b| < a $, $ \kappa  > 1$; and $\kappa = 1$ for rolling over a plane ($a \rightarrow \infty$).
One can also consider the case $- \infty < \kappa < 0$, when the fixed ball
$S_a$ lives inside the rolling ball $S_{|b|},  b < 0,  |b| > a $.)

The formula for  the energy is the same as in the sphere-plane case,
\begin{equation}  \label{totalenergy}
  2H = 2T =  \mu b^2 || \Omega \times \gamma ||^2 +   I_1 \Omega_1^2 + I_2 \Omega_2^2 + I_3 \Omega_3^2  \,\,.
 \end{equation}
We now give an argument to explain why the invariant measure stays the same,
independently of the value of $\kappa$.
\begin{proposition}
There is  a smooth invariant measure for the marble rolling of a sphere over
a sphere.   The density function is  the same as in the case of Chaplygin's sphere rolling over a plane, $F(\gamma)^{-1/2}$, with
$F$  given by (\ref{chapfunction}).
\end{proposition}
\noindent {\bf Proof.} As a warm up, let us verify first that the  total energy function 
(\ref{totalenergy}) of the sphere rolling
 over a sphere (a formula in which 
the radius of the base sphere does not appear)    is conserved along the trajectories of (\ref{L},
\ref{kapa}).  
This  direct check will be instructive, since all sphere-sphere systems have the same
energy function and the same $ \dot{L} $ equation, but a {\it family } of $\dot{\gamma}$ equations,
parametrized by $\kappa$. We still get
$ \dot{H} = {\rm grad}_{(L,\gamma)} \, H \cdot (\dot{L}, \dot{\gamma}) \equiv 0$ because  {\it both}
\begin{equation} \label{nablaL}  \nabla_L \, H \cdot \dot{L} \equiv 0  \,\,,
\end{equation}
\begin{equation}  \label{nablagama} \nabla_{\gamma} \, H \cdot \dot{\gamma} \equiv 0   \,\,.
\end{equation} 
 Equation (\ref{nablaL}) is immediate, since $\nabla_L \, H = \Omega $
and $\dot{L} = \Omega \times L$ is perpendicular to $\Omega$.  Now, we know that the proposition  
holds for the case $\kappa = 1$.  Equation (\ref{nablagama}) must be therefore true for $\kappa = 1$. 
Now, for an arbitrary $\kappa$,  the left hand side of (\ref{nablagama}) only picks this multiplicative factor
so it {\it must } also be  zero and there  is no  need to do the computation).  

In order to verify that the measure is invariant, we follow a similar procedure.
Let us scrutinize the 
derivation presented Duistermaat,  \cite{Duis} (lemma 7.1) for the case $\kappa = 1$.  His approach was to split
$$  {\rm div}_{(\gamma,L)}\, ( \gamma \times \Omega(\gamma,L)\,, L \times \Omega(\gamma,L) )
$$ 
as the trace of  the derivative with respect to $\gamma$ of the first component
plus the trace of  the derivative with respect to $L$ of the second component.
The latter trace is obviously identically zero, whereas the
first one is found to be (see (7.2) in that paper)
\begin{equation}  \label{eqdiv}
{\rm div}_{\gamma} (\gamma \times \Omega(\gamma,L))) =  \frac{(\gamma \times \tilde{A}^{-1} \gamma, \tilde{A}^{-1} L)}{F(\gamma)} \,\,.
\end{equation}
 On the other hand,  the directional derivative of $F(\gamma)$  is found to be
 \begin{equation}  \label{eqder}
  ({\rm grad} F , \gamma \times \Omega(\gamma,L)) = 2 (\gamma \times \tilde{A}^{-1} \gamma, \tilde{A}^{-1} L)
\end{equation}
and hence
 $$ {\rm div}_{(\gamma,L)} \, [ F^{-1/2} ( \gamma \times \Omega(\gamma,L)\,, L \times \Omega(\gamma,L) )] = {\rm div}_{\gamma} \, [ F^{-1/2}  \gamma \times \Omega(\gamma,L) =  0   \,\,.
 $$
{\it  Now, comes the argument: the only change caused by the presence of $\kappa$ is to multiply both (\ref{eqdiv}) and (\ref{eqder}) by this same fator, so the net result stays zero.}
 
  Functions $f_1, f_2, f_3$ are still  integrals of motion for any $\kappa$, but the ``area integral''  $f_4$ is lost! It is easy to see that it holds 
 only for $\kappa = 1$.  In the three dimensional manifolds $M^3$ defined by level sets of the three (surviving) integrals we have an 
 invariant measure.  Flows on  three dimensional compact manifolds with an invariant measure are of interest for  dynamicists
and one  expect that
the motion will be ergodic in $M^3$ except perhaps for some special values of $k$,  where some new integral of motion may  appear. 

Borisov found the only known (perhaps unique) new integrable case,  $\kappa = -1$, corresponding to an internal fixed ball of half
the radius of the rolling ball,
see  \cite{BM}, pag. 194.

\begin{remark}  The external case with equal radii (corresponding  to $ \kappa = 1/2 $) could be (erroneously) thought to be also integrable,  for a simple reason.  As we will see below, the two dimensional {\it rubber} distribution is {\it
holonomic} for two spheres of equal radius. {\it One could hope that relaxing the no-twist 
constraint would keep the distribution holonomic}, as it amounts to allow any
finite rotation of the moving sphere about the normal at any given (but arbitrary) point on the base sphere. 
Each leaf ${\cal F}$ would then be three dimensional,
 forming a
a  $S^1$ bundle over $S^2$, and the topology would be (probably) that of $SU(2)$.
Unfortunately there is a serious blunder is this reasoning. Although ${\rm span} (X_1, X_2)$
is integrable, adding the infinitesimal
rotations about the contact points produces a 3 dimensional distribution
${\rm span} (X_1, X_2, X_3)$  with growth  3-5.
\end{remark}

\section{Rubber rolling over a sphere as a  $SO(3)$ Chaplygin system.}

\subsection{ Reduction to $T^*S^2$.}    

Consider the principal bundle  
with total space  $Q = S_a \times SO(3)$,  
$$
 SO(3) \hookrightarrow   Q  \rightarrow  Q/SO(3) \,\,\,,\,\,\,  R \cdot  (q, S) = (R(q), RS ) 
$$
given  by the diagonal action.  It is not hard to see that the no-slip and no-twist constraints
define a principal bundle connection for any rolling surface $\Sigma_2$. We call it the {\it rubber connection}.

Recall (\ref{connec1}),  the Ehresmann connection for the rubber rolling of a  surface $\Sigma_2$   on an arbitrary base $\Sigma_1$:
$$
\omega = N_1(q_1) \times ([ dN_{g \Sigma_2} + dN_1 ](q_1) \cdot \dot{q}_1   \,\,.
$$

When $\Sigma_1$ is a sphere there is an extra feature, namely, {\it equivariance}, so we have a principal bundle connection.  We need to show that  if $\omega$ solves both 
equations  (\ref{notwist}, \ref{rolling2}), which in this case become
$$  \omega \times q_1 = [I + a dN_{g \Sigma_2} (q_1)] \dot{q}_1\,\,\,,\,\,\, \omega \cdot q_1  = 0 \,\,\,,
$$ 
then  $R \omega$ solves
$$ R \omega \times R q_1 = [ I + a dN_{Rg \Sigma_2} (R q_1)]( R  \dot{q}_1)  \,\,\,\,,\,\,\,\, \R \omega \cdot R q_1 = 0 \,\,.
$$
The second fact is obvious, and  the first  follows from\footnote{Latexing  is more difficult 
than  the result. A drawing should make it obvious.}
$$   dN_{Rg \Sigma_2}(Rq_1) =  R  dN_{g \Sigma_2} (q_1) R^{-1} \,\,\, . $$


See   \cite{KoillerArma} for  reduction   written in  Lagrangian form. The reduction procedure in
Hamiltonian form, following the recipe  given  in \cite{Kurt2} is as follows: 
\begin{theorem}   \label{reductionchaplygin}
Let $\Sigma_2 = \Sigma$ be arbitrary but  $\Sigma_1$   a sphere of radius $a$. Then rubber rolling is described by  a non-abelian  Chaplygin system on the bundle $Q = SO(3) \times S^2(a)$, with
 symmetry  group $SO(3)$ acting diagonally.  The constraints are given by 
 the rubber connection
\begin{equation}  \label{connectiongeneral}
\omega = \frac{q_1}{a} \times [ dN_{g \Sigma_2} \cdot \dot{q}_1 + \frac{\dot{q}_1}{a} ]  \,\,.
\end{equation}
Define the projection is
 by  $ \pi (R, q) = R^{-1} q/a  \in S^2  \,,\, R \in SO(3)\,,\, q \in S^2(a)$.
Then the dynamics reduces to an almost-symplectic system  $(H,\omega_{NH})$ in  $T^*S^2$. 
The Hamiltonian $H$  in $T^*S^2$ is the Legendre transform of the  compressed Lagrangian $T(\gamma,\dot{\gamma})$ given by
\begin{equation}
2T = (\tilde{A} \Omega,\Omega)\,\,\,{\rm with}\,\,\,\,\,\,\, \Omega = [ ( Id + \frac{1}{a} \,(dN_2)^{-1}(\gamma) )] (\dot{\gamma})  \times \gamma \,\,\,\,,\,\,\,
\tilde{A} = A + \mu ||Q||^2 id \,\,\,,\,\,\, Q = - N_{\Sigma}^{-1} (\gamma)
\end{equation}
obtained from the original one by horizontal lifting via the connection. 
The  2-form  is given by  $ \omega_{nh} = \omega^{T^* S^2}_{can} + (J,K) \,\, $ and in general is non closed.  
The ammended  term $ (J,K)$  is semibasic;  $J$ is the momentum  (\ref{momentummap})
of the $ SO(3)$ action and  $K$ is the curvature (\ref{Ehrcurv}) of the rubber connection. 

 The domains of $J$ and $K$ are matched by the Legendre transformation from $TQ$ to $T^*Q$.  Their  $ Ad$  and 
$Ad^*$ ambiguities cancel out automatically. In some more detail, 
\begin{equation}
(J,K) = ( 1 - \frac{1/a^2}{k_{\Sigma}(\gamma)} ) \,\, (\gamma, M ) \, d \,{\rm area}_{S^2} \,\,,
\end{equation} 
where  $\gamma = q/a$ is the contact point in the base sphere,  $k_{\Sigma}(\gamma)$
is the Gauss curvature of the rolling surface at the current contact point
$Q = - N_{\Sigma}^{-1} (\gamma)$,
and $M = A\Omega$ is the angular momentum at the center of mass of the rolling body.
\end{theorem}

We plan to use this result as a basis for future work.  The first in line is the rubber rolling over a sphere of
surfaces of revolution with two equal inertias.
For the rest of this paper we restrict ourselves to the sphere-sphere example.

 \subsection{The rubber sphere-sphere connection 1-form. } 
 
  We show that for the sphere-sphere problem there is
an invariant measure, and even better, the reduced system is conformally symplectic\footnote{We conjecture that this is also true for bodies of revolution.}.
The expert can skip many of the redundant calculations below, intended as checks of the general
results presented earlier. When   $\Sigma_2$ is a sphere of radius $b$,  (\ref{connectiongeneral}) 
 becomes
\begin{equation}  \label{connec3}
 \omega = \frac{1}{a^2}\, (1 \pm  \frac{a}{b}) \, q_1 \times  \dot{q}_1  \,\,.
 \end{equation}
 Recall that for simplicity we drop the $\pm$ sign, using the convention that $b$ is negative one of the spheres contains the other.
\begin{lemma} \label{horizontallift}
The horizontal lift of  $\,\, \dot{\gamma} \in T_{\gamma} S^2 \,\, $ at $\,\, (a R(\gamma), R) \, \in Q \,\,$ is given by
\begin{equation}   \label{Hor}
\dot{\gamma} \mapsto Hor(\dot{\gamma}) = ( - b R(\dot{\gamma}),\, (1 +   \frac{b}{a})\,  (L_R)_* [\dot{\gamma} \times \gamma] )  \,\,.
\end{equation}
where $(L_R)_*$  indicates the left translation and $[ , ]$ the standard isomorphism
from $\R^3$ to $sO(3)$.
\end{lemma}
\noindent {\bf Proof.}  The formula looks odd  due
to the fact that the $SO(3)$ action is diagonal.   Note first that
$$ (\pi_{*})_{q,R)} \, (\dot{q}, \dot{R}) = \frac{1}{a} ( R^{-1} \dot{q} - R^{-1} \dot {R} R^{-1} q ) \,\,\,.
$$
Take  $ a R(\gamma) = q, \dot{q} = - b R(\dot{\gamma})$  and  $ \dot{R} = (L_R)_* [\dot{\gamma} \times \gamma] $.  Then
$$ (\pi_{*})_{(q,R)} \, (\dot{q}, \dot{R}) = - \frac{b}{a} \dot{\gamma} -
(1 +   \frac{b}{a}) (\dot{\gamma} \times \gamma) \times \gamma ) = \dot{\gamma}
$$
and a short calculation shows that (\ref{connec3}) holds:
$$ \frac{1}{a^2}\, (1 + \frac{a}{b}) \,  R^{-1}  (q_1 \times  \dot{q}_1 ) =
 \frac{1}{a} (1 + \frac{a}{b}) \gamma \times  (- b \dot{\gamma} ) = 
( 1 + \frac{b}{a}) \dot{\gamma} \times \gamma \,\,.
$$
In order to find the connection form, it sufficies to compute along the global section  $\gamma \mapsto (a \gamma, I)$.    

\begin{proposition}  The connection 1- form $\phi$ interpreted as a 1-form in $S^2 \times SO(3)$ with values in $s0(3) \equiv \R^3$,  is given  at  the point $ ( \gamma, I)$ of the global section  by
\begin{equation}  \label{connectionform}
  \phi (\dot{\gamma}, \sigma) = (1 + b/a) \gamma \times \dot{\gamma} - (b/a) [ \sigma - (1 + a/b)\, (\sigma \cdot \gamma)\, \gamma ]   \,\,.
\end{equation}
\end{proposition}
{\bf Proof.}
 By a direct check one verifies, as desired, that
$$ \phi (\rho \times \gamma, \rho ) = \rho  \,\,\,\, , \,\,\,\,
 \phi ( - \frac{b}{a}) \dot{\gamma}, ( 1 + \frac{b}{a}) \, \dot{\gamma} \times  \gamma) = 0 \,\,.
 $$
 
While the metric is {\it left } invariant under the diagonal action of $G = SO(3)$ on
$S_a \times SO(3)$, looking at (\ref{connec3}) we observe that 
in the sphere-sphere problem the constraints are  invariant under 
another action of $SO(3)$ on $S^2 \times SO(3)$, namely the {\it right} action on the group component:   $ (\gamma , R) S = (\gamma, RS)$.
This makes the problem reminiscent to a  LR nonholonomic Chaplygin  system
$ H \hookrightarrow G \rightarrow H/G $,  see \cite{Fedorov}, except that here
the bundle structure is of the form
\begin{equation}  \label{bundle}
G \hookrightarrow G \times H/G \rightarrow H/G 
\end{equation}
with $H \equiv S^1$ the rotations about the vertical axis, $H/G \equiv S^2$ the
homogeneous coset space.

 Assuming that
   results similar to those for LR systems would still hold,  since the reduced  dynamics has two degrees of
freedom,
we predict the existence of a function $f: S^2 \rightarrow \R$ such that  $ d(f \omega_{nh}) = 0 $. 
In other words, the reduced
system in $T^*S^2$ would be Hamiltonizable,  since the nonholonomic 2-form is conformally symplectic\footnote{See \cite{Haller,Haller1,Vaisman,AissaWade}
to appreciate the consequences of being conformally symplectic.} and in particular,  there will be a smooth  invariant measure.
Actually, we guess $f$ from the  density function  of the invariant measure,  relative to the
Liouville volume of $T^*S^2$.  The good news is that it is possible to verify directly if a nonholonomic system
has a smooth invariant measure, see \cite{Blackall,Cantrijn,KobayashiOliva,KupkaOliva,Blochnonlinearity}.




 In future work we plan to pursue a LR theory for principal  bundles of the type (\ref{bundle}). Here, we use an alternative, albeit tentative approach: in order to make the LR analogy more plausible, we consider   the whole  ``onion'' of base
 spheres $S_a, a > 0$ 
(the singular  limit $a = 0$  represents Veselova system).  Then our  mechanical  system
lives on an {\it extended} configuration space,  $ Q  =  E(3)^{\bullet} = SO(3) \times (\R^3 - 0) $  with  a metric which is
 manifestly    invariant under the {\it  left}  action of $SO(3)$ :
$$
  2T = \mu (1 \pm b/|q|)^2 || \dot{q} ||^2  +   I_1 \Omega_1^2 + I_2 \Omega_2^2 + I_3 \Omega_3^2
$$
Note that the presence of $|q|$ in the denominator forces us
to remove the origin  of $\R^3$.

The  pseudo LR system has bundle structure\footnote{Using a more abstract language helps to generalize to configuration spaces that are semi-direct products of Lie groups.}
\begin{equation}  \label{bundle2}
   SO(3)) \,\,  \hookrightarrow \,\,    E(3)) \,\,\,  \rightarrow   \R^3  \,\,\,,\,\,\,  R \cdot  (S,q) = (RS, Rq ) 
\end{equation}

We now rewrite the no-slip, no-twist  conditions in terms of differential forms. The right  and left invariant forms in $SO(3)$ will be  denoted here
\begin{equation}
  dR R^{-1} = \left( \begin{array}{lll}   0 & - \rho_3  &  \rho_2 \\  *  & 0 &  - \rho_1 \\ *  &  *  & 0 
\end{array} \right)  \,\,\,,    \,\,\,   R^{-1} dR =  \left( \begin{array}{lll} 0 & - \lambda_3  &  \lambda_2 \\  *  & 0 &  - \lambda_1 \\ *  &  *  & 0 
\end{array} \right) \,\,.
\end{equation}

The distribution
${\cal D}$ of 2-dimensional subspaces  in $E(3)^{\bullet}$
 with  $ |q| \neq 0$,  is the joint  kernel of the four differential forms $\eta_i, \, 1 \leq i \leq 4$  in  $ E(3)$ defined by 
\begin{equation}
(\eta_1, \eta_2, \eta_3)^t =  \vec{\rho} - \frac{1}{|q|^2} \, (1 + \frac{|q|}{b} ) q \times dq \,\,\,,\,\,\,\,  \eta_4 =  q_1 dq_1 + q_2 dq_2 + q_3 dq_3  \,\,.
\end{equation}
where   $\vec{\rho} = (\rho_1, \rho_2, \rho_3)^t $.  Notice that the Pffafian equation  $(q,\rho) = 0 $  defining the no-twist condition is  satisfied 
automatically when $\eta_i=0, 1 \leq i \leq 4$. 
Notice also that $S \in SO(3)$ acting on the right on an element $(R,q) \in SE(3)$ gives $(RS,q)$.
This is best seen using the costumary representation
for $E(3)$:
$$  (R, q)  \leftrightarrow  \left( \begin{array}{ll}  R & q \\ 0 & 1
 \end{array}
   \right)
$$
This implies 
\begin{proposition} Consider the distribution $ {\cal D}$ in $E(3)^{\bullet}$ defined by the constraints  $\eta_i = 0\,,\,\,  \, 1 \leq i \leq 4$.  These forms   are invariant under the right action  
of $SO(3)$ in $SE(3)^{\bullet}$, that is $  {\cal D}_{(R,q)} =  {\cal D}_{(I,q)} \cdot R $.
\end{proposition}

\subsection{Curvature of the sphere-sphere connection}  \label{curvsphsph}

This section could be avoided in view of (\ref{Ehrcurv}) but we decided to
include it for completeness. We already know the connection 1-form $\phi$ at the global section. In order to compute the curvature 2-form we need  an expression everywhere. This can be done by equivariance,  but this is a somewhat unpleasant task [ $\phi_{(q,R)}(\dot{q},\dot{R})$ is obtained from (\ref{connectionform})
by computing  $R$ of the result of $\phi  (\dot{\gamma}, \sigma) $ with the inputs  $\gamma \leftarrow  R^{-1}q/a, \, \dot{\gamma} \leftarrow  R^{-1} \dot{q}/a - \omega \times \gamma, \,  [\sigma ]  \leftarrow  R^{-1} \dot{R} \,,\,\, \omega = \dot{R} R^{-1}$.
Also, we must extend it to a connection in $\R^{\bullet} \times SO(3) \hookrightarrow SE(3)^{\bullet} \rightarrow S^2$, by adding  a component $(\dot{q}, q/|q|)$
to $\phi$, and  replacing  $a$ by $ |q|$  and $ \dot{q}$ by
$ \dot{q} - (1/|q|^2) (q,\dot{q})$.  ]
  
  After putting everything inside $\phi$, we would get an unpleasant formula. Computing the curvature
directly would be a ``tour de force''. We prefer to guess the curvature by
a "wrong" but algebraically elenant calculation. We then reinterpret and confirm the guess
using elementary geometry.

 We can use cartesian coordinates in $\R^3$, since we extended the base space to   an ``onion'' of spheres.
Consider the vectorfields in $ SE(3)^{\bullet} $ given by
$$ Y_1 = \partial/\partial x +  (\frac{1}{a}  + \frac{1}{b}) (\frac{z}{a}\,  X_2 - \frac{y}{a} \, X_3)\,\,, \,\, Y_2 = \partial/\partial y +    (\frac{1}{a}  + \frac{1}{b})   (-  \frac{z}{a}\, X_1 + \frac{x}{a} \,  X_3)$$
where $X_1,X_2,X_3$ are the right invariant vectorfields of $SO(3)$ ($[X_1,X_2] = - X_3$, etc). 
 These  vectors are horizontal at the point $ q_a  = (0,0,a),\, R =  I $.  A simple mental calculation gives  us
$$ [Y_1,Y_2]_{|(q_a,I)} =  [2 \frac{1}{a}  (\frac{1}{a}  + \frac{1}{b})  -  (\frac{1}{a}  + \frac{1}{b})^2  ] \,X_3  \,\,\, ,
$$
$$  [2 \frac{1}{a}  (\frac{1}{a}  + \frac{1}{b})  -  (\frac{1}{a}  + \frac{1}{b})^2  ] =  \frac{1}{a^2} - \frac{1}{b^2}  \,\,\,.
$$ 
Now, both $Y_1$ and  $Y_2$ must be multiplied  by  $-b$ in order to be the
horizontal lifts at $q_a$ of the unit vectors  $(1,0,0)$ and $(0,0,1)$
at  the north pole of  $S^2 = S^2(1)$, see  equation (\ref{Hor}) in lemma \ref{horizontallift}.  As a result, we have to multiply that result by  $b^2$.
From H.Cartan's magic formula, we get (note the extra - sign):
\begin{proposition} The  curvature of the connection, interpreted as a 2 form in the unit sphere $S^2$ with values in $\R^3 \equiv sO(3)$, is
\begin{equation}  \label{curvatsphsph}
 K(\gamma) = ( 1 - \frac{b^2}{a^2} ) \, \gamma  \, \,d {\rm area}_{S^2} 
\end{equation}
\end{proposition}
This calculation  can  be made more rigorous.  We forgot to compute   $ Y_1 \phi(Y_2) - Y_2  \phi(Y_1)$  in H. Cartan's  ``magic formula'',
 at $(0,0,a)$.   By a symmetry argument,  one should be able to show  that it must
 vanish.

We now present a geometric way to compute the curvature.  Consider a  small circle $C_1$ in the base sphere $S_a$.  Roll  along $C_1$, without slip nor twist, a ball
$S_b$ (for concreteness, we do it externally),  until it  comes back to the initial point in $C_1$. The curvature is the following limit:
\begin{center}
$\lim \,\,\,$ ( how much a frame attached to $S_b$
has rotated )  $\,\,\, /   \,\,\,\,$  ( area of the cap in $S_a$  )  $\,\,\,\,\,\,\,\,\,\,, \,\,\,\,\,$  as  $S_a$  shrinks . 
\end{center}
Now,  Proposition \ref{geosame} tells that the corresponding curves have the same geodesic curvatures, so the  contact point on $S_b$ must move in another circle $C_2$. After one turn arond $C_1$, the end point in $C_2$ is not the
same as the initial point. For instance, the marker does  several turns plus a fraction if $b$ is smaller than $a$.   Moreover,
proposition \ref{propstruik} allow us to replace the spheres by their 
enveloping cones\footnote{Cone rolling  has many applications, ranging from transmission gears to  leisure stuff, 
see http://www.kurt.com/gear.html\,\,, (no relation with  one of the authors), http://www.1percent.com/store/cart/RLZENCONE.html 
(no relation with any of the authors,
and   http://kmoddl.library.cornell.edu/model.php?m=reuleaux,(yes, we enjoy anything  related to Cartan's 1:3 rolling! }.

 The equivalent problem is the following: let a cone $\Sigma_2 $   roll  around   another cone  $\Sigma_1$   until its generator   returns to the same place.  How much 
  the position of a marker point
 $Q_2$ in $\Sigma_2$ has  rotated around axis $\ell_2$ ?  We call this angle
 divided by $2 \pi$ the {\it holonomy}.
  The reader should be able to sketch  a figure. The moving cone $\Sigma_2$ rolls around
the fixed cone $\Sigma_1$,  having a common a different generatrix {\cal g} at every instant.   
 Draw right triangles $V O_1 q$ and $V O_2 q$ where   $V$ is the common vertex,  $O_i$ in the axis of $\Sigma_i$ are the centers of the circles $C_i$, and
 $q$ is the contact point. Let $\alpha = \angle O_1 V q$ and $\beta = \angle O_2 V q$ the cones apertures, $ \theta = \pi/2 - \alpha = \angle O_1 q V $,
 $ \phi = \pi/2 - \beta = \angle O_q q V $.



  The answer is  easy:  
\begin{center} $ |{\rm holonomy}| =  \,\,\, $   length of the base circle $C_1$  $ \,\,\,\,  /   \,\,\,\, $   length
of $C_2  \,\,\,= \sin(\alpha)/\sin(\beta) = \cos(\theta)/\cos(\phi) \,\, .$
\end{center}

Now by elementary trigonometry,
\begin{equation} \label{thetaphi}   a \tan(\theta) = b \tan(\phi)  \,\, ( = \ell) \,\,.
\end{equation}
which by the way says that the geodesic curvatures are
 \begin{equation}   \kappa_{C_1} =   \kappa_{C_2} =  1/\ell  \,\, .
\end{equation}
   Thus
\begin{equation}   {\rm holonomy}  =  -  \frac{a}{b} \, \frac{\sin(\theta)}{\sin(\phi)}  \,\,\, , 
\end{equation}
 Moreover, from  (\ref{thetaphi})  we can  solve  for  
$$\sin(\phi) = \frac{a}{b} \frac{\tan(\theta)}{\sqrt{ 1 + \frac{a^2}{b^2} \tan^2 (\theta)}}  $$
so that  
\begin{equation}   {\rm holonomy} = -  \cos(\theta) \sqrt{ 1 + \frac{a^2}{b^2} \tan^2(\theta)}  \,\,.
\end{equation}
Why the  minus sign?  The angle between the normals is $\alpha + \beta$, and since  $\theta + \alpha = \pi/2,\, \phi + \beta = \pi/2$,  both $\alpha$ and $\beta$ tend
to  $\pi/2$ when  $\theta$  (and so $\phi$)  tends to zero.  So the two normals tend to become {\it opposite} to each other !  
Actually, here is one more twist:   a frame  attached to the  sphere $S_b$ rotates by    $2 \pi - 2 \pi |{\rm holonomy}| $.
The extra $2 \pi$ is due to the ``Phileas Fogg effect'' (known to mathematicians as Hopf's  {\it umlaufsatz}).  Now, the area of the cap in sphere $S_a$ is
$ 2 \pi a^2 ( 1 - \cos(\theta) ) \sim  \pi a^2 \theta^2 $.

In short, we end up needing to compute the limit
$$  \lim_{\theta \rightarrow 0} \,\, \frac{2}{a^2} \,  \frac{1 - \sqrt{ 1 + \frac{a^2}{b^2} \tan^2 (\theta)} \, \cos(\theta)}{\theta^2}
$$ 
This is a nice calculus exercise, the result being
that  the curvature is  inded  $\frac{1}{a^2} - \frac{1}{b^2}$ times the infinitesimal rotation about the normal around the contact point $q$ in sphere $S_a$. 

\begin{remark}  Among the ratios $a/b$, two are exceptional: $a=b$ provides
a holonomic distribution since the curvature vanishes.  Cartan has used
the 3:1 ratio to construct his exceptional Lie group $G_2$, see \cite{Cartan,MontBor}.
We hoped that the corresponding rolling problems would have some special
features, but that was with no avail. 
\end{remark}


 \subsection{Hamiltonization of the  sphere-sphere system } 
  \label{Hamiltonization}

We will compute step by step  
 the 2-form $\omega_{NH} =  d p_{\gamma} d \gamma + (J,K)$, and the reduced energy function
$ H: T^*S^2 \rightarrow \R$,  to convey   the algorithmical 
nature of the method.  
\begin{itemize}
\item  The momentum mapping is given by  $ J(m_q,r) =  m + q \times r  \,\,$, see equation (\ref{momentummap}).
\item  Pairing the momentum with the curvature (\ref{curvatsphsph})   along the global section $ q \in S^2(a) \mapsto (q, I), \,\, q =  a \gamma $ we get
$$   (J,K) =   ( 1 - \frac{b^2}{a^2} ) \, ( \gamma,   M ) \, d {\rm area}_{S^2} \,\,\,,\,\,\, M = A \Omega
$$ 
\end{itemize}

In order to match the domains, we   follow the ``clockwise diagram'' presented in \cite{Kurt2}. This diagram is a consequence of a Hamiltonian counterpart of
Lagrangian reduction \cite{CeMaRa}
using quasi-coordinates \cite{KoillerRomp}.
\begin{itemize}
\item  The Legendre transform $TQ \mapsto T^*Q$  induced by the Lagrangian $ T =  \frac{1}{2} [  \mu  (1 + \frac{b}{a})^2\, |\dot{q}_1|^2 + (A \Omega, \Omega) ] $, at the point
$  (\dot{q}_1, \dot{R})_{(q_1, R)} = (v , \Omega)_{(a \gamma, I)} \in  T_{(a \gamma, I)} Q$ is $(P_v, P_I)_{(q_1, R)}$, given by 
$$ P_I  = M  =  A \Omega \,\,\,, \,\,\,\, P_v = r =  \mu \, (1 + \frac{b}{a})^2\, v   \,\,\, .
$$
\item   The Legendre transform of the horizontal lift
 $Hor(\dot{\gamma})  =  ( v, \Omega ) =  ( - b \dot{\gamma}, (1 + \frac{b}{a})  \dot{\gamma} \times \gamma \,) $:
$$  m  = (1 + \frac{b}{a}) A ( \dot{\gamma} \times \gamma) \,\,\, , \,\,\,\, r =  \mu \, (- b)(1 + \frac{b}{a})^2   \dot{\gamma} \,\,\,\,
$$
\end{itemize} 
The end result is
\begin{equation}  \label{JKab}
 (J,K) = (1 - \frac{b^2}{a^2})\,  (1+\frac{b}{a}) (\gamma, A( \dot{\gamma} \times  \gamma) \,  dS^2  =  (1 - \frac{b}{a}) (1+\frac{b}{a})^2 (\gamma, A( \dot{\gamma} \times  \gamma) \,  dS^2 \,\,.
\end{equation}
\begin{itemize}
\item To obtain the compressed Lagrangian  $ T_{red}$, we compute the Lagrangian at $Hor(\dot{\gamma})$:
\begin{equation}
T_{red} = \frac{1}{2} (1 + \frac{b}{a})^2 [  \, \mu \,b^2 \dot{\gamma}^2 +  ( A(\gamma \times \dot{\gamma}), \gamma \times \dot{\gamma} )\, ] 
\end{equation}
\item
 The Legendre transform in the base,   $\,\,\, p_{\gamma} = \partial T_{red} / \partial \dot{\gamma} = G(\gamma) \cdot \dot{\gamma} $:
\begin{equation}   \label{legendreab}
p_{\gamma} = G(\gamma) \cdot \dot{\gamma}  = (1 + \frac{b}{a})^2 [ \mu b^2  \dot{\gamma} +      A(\gamma \times \dot{\gamma)}\times \gamma  ]
\end{equation}
\item  The determinant of $G(\gamma)$, restricted to $T_{\gamma} S^2$ is
\begin{equation}
F(\gamma) = [{\rm det} \, Leg_{comp}] = (I_1I_2I_3)(1 + \frac{b}{a})^2 \, \left( (A^{-1} \gamma, \gamma) + 
\mu\,b^2 [ \frac{\gamma_2^2 + \gamma_3^2}{I_2I_3} + \frac{\gamma_1^2 + \gamma_3^2}{I_1I_3} + \frac{\gamma_1^2 + \gamma_2^2}{I_1I_2}] + \frac{\mu^2 b^4}{I_1I_2I_3} \right) \,\,.
\end{equation} 
\end{itemize}
Only this last assertion represents a nontrivial calculation, but that was 
done in \cite{Kurt2},  for the limit  $a \rightarrow \infty$ of the rubber ball rolling over a plane. 
In that paper, we announced ((3.38), Proposition 3.6 of \cite{Kurt2}):
\begin{proposition}   For  $ \, a = \infty \, $ (rubber ball rolling over a plane)
\begin{equation}
 d (  f ( \omega_{can}^{T^*S^2} + (J,K) \, ) = 0  \,\,,\,\,\, f = F^{-1/2} .
\end{equation}
\end{proposition}
\noindent {\bf Proof.} This was  checked directly by computer algebra, using spherical coordinates.  It was useful to pull back the
canonical two form in $T^*S^2$ to $T S^2$ via  $p_{\gamma} = G(\gamma) \cdot \dot{\gamma}$.   
We can provide a  Mathematica notebook for those interested.

We now show how to adapt this result  for  $b/a  \neq 0$.
First note that the conformality condition 
 can we rewriten, {\it when the base space of the bundle  has dimension two},  as
$$  d (\log f) \wedge \omega_{can} + d(J,K) = 0   \,\,.
$$
From  (\ref{JKab}) and (\ref{legendreab})  we observe that $\omega_{can}$, when pulled back to $TS^2$, picks the  factor  $ (1 + b/a)^2$, whereas $(J,K)$ picks
$(1 - b/a) (1 + b/a)^2$, that is,  one fator $1 - b/a$ surplus.   Thus we can't have the same $f$ as in the
case $a = \infty$, but   if we change  $f$ to  $f^{1- b/a}$, that should work.
\begin{theorem}   The reduced rubber ball rolling   over a ball is
conformally symplectic, 
\begin{equation}  \label{symplecticform}
 d  \, [  f_{a,b}(\gamma) ( \omega_{can}^{T^*S^2} + (J,K)  )\, ] = 0   
\end{equation}
with conformal factor  $$ f_{a,b}(\gamma) = F^{\frac{b-a}{2a}} \,\, .$$
\end{theorem}
\noindent {\bf Proof.} For the suspicious reader, we can provide a Mathematica notebook  using spherical coordinates.  But  a
 derivation (done by hand) follows from sphero-conical coordinates, as we show
 in the
next section.
\begin{remark}
  The presence of the exponent  $\frac{b-a}{2a}$ is somewhat mysterious, but it makes sense: 
when
$b = a$,  then  $ f \equiv 1 $ (as it should be) since in that case
$(J,K) \equiv 0 $ and the problem becomes Hamiltonian.
The predicted exponent from the results   
 by  Veselov and Veselova  \cite{Veselovs} or 
Fedorov and Jovanovic \cite{Fedorov} would be  $-1/2$, but this
corresponds only to  $b/a \rightarrow 0$.   
We believe that the
problem resides on the singularity at the origin $a=0$, and the mystery will
be resolved with a thorough study of bundles (\ref{bundle}).  
  
\end{remark}



All these problems have symplectic forms  parametrized by $(a,b)$ and the same  Lagrangian on $TS^2$ containing two terms, 
of the form $ \dot{\gamma}^2 $ and $( A(\gamma \times \dot{\gamma}), \gamma \times \dot{\gamma} )$, namely 
\begin{equation}  \label{lagrangian rubber}
  2T/(1 + b/a)^2  =  \mu b^2 \dot{\gamma}^2 +  ( A(\gamma \times \dot{\gamma}), \gamma \times \dot{\gamma} )  \,\,.
\end{equation}


For $a/b = 1$ (the case of equal spheres)  we have  $f \equiv 1$ and $(J,K)=0$, so we have the standard
symplectic form in the base $T^* S^2$.
Based on the classic integration by Jacobi of the ellipsoid geodesics problem,
we could hope that the the Hamiltonian   separates  in sphero-conical coordinates $\lambda_1,\lambda_2$, so the Hamiton-Jacobi
method would work.  We were be lead to a disgusting failure.

\subsection{Sphero-conical coordinates}

This section is also standard material for the experts. Sphero-conical coordinates appear naturally 
if one wants to express (\ref{lagrangian rubber}) in the form
$
2T = c_1 \dot{\lambda}_1^2 + c_2 \dot{\lambda}_2^2  \,\,.
$
To achieve this aim we need   to diagonalize  the positive quadratic form  $ (Ax,x)\,,\,\,\, x \in \R^3, \,\,A = diag(I_1,I_2,I_3)$  
 when restricted to  subspaces defined by  $(x,\gamma) = 0$.  
This is a well known problem, 
going back to Lord Rayleigh.  It boils down to
extremize $(Ax,x)$ subject to $\sum x_j^2 = 1$ and $(x,\gamma) = 0$. Let us apply the
Lagrange multiplier method, with an ammended objective function 
\begin{equation} \label{U}
 U = \frac{1}{2} \, \sum \, I_j x_j^2 - \lambda (\sum x_j^2 - 1) - \mu \sum x_j \gamma_j  \,\,.
\end{equation}
The first order condition $  \nabla U_x = 0 $  yields
\begin{equation}  \label{xj}
 x_j = \mu \frac{\gamma_i}{I_j - \lambda}, \,\,  1 \leq j \leq 3 .
\end{equation}
Inserting in the constraints $\sum x_i^2 = 1 \,,\,\, \sum \, x_i \gamma_i = 0  $
we get 
\begin{equation}     \label{conical}
\sum \, \frac{\gamma_j^2}{I_j - \lambda} = 0
\end{equation}  \label{muvalue}
\begin{equation}  \mu^2 = 1/[ \sum \frac{\gamma_j^2}{(I_j - \lambda)^2} ]  \,\,.
\end{equation}
This holds in $n$ dimensions, but let us keep working with $n=3$ for concreteness. 
If we order the parameters as  $I_1 < I_2 < I_3$, equation (\ref{conical}) has two roots,
$$  I_1 < \, \lambda_1 = \lambda_1(\gamma;I_1,I_2,I_3) \, < I_2 $$   
$$  I_2 < \, \lambda_2 = \lambda_2(\gamma;I_1,I_2,I_3) \, < I_3 $$ 
(they can be found explicitly as the roots of   a second degree equation in $\lambda$ thinking of the $\gamma_j$ as parameters).  

 For each choice of $\lambda$, (\ref{conical}) represents an elliptical
cone in the three dimensional space $(\gamma_1,\gamma_2,\gamma_3)$. It is easy to see that the curves of intersection of these cones with 
the sphere  $ (\gamma,\gamma) = 1 $ form an orthogonal
system of coordinates.  This is because we know that the two vectors (\ref{xj}) that solve (\ref{U}) are perpendicular, and moreover, these vectors
are manifestly parallel to the gradients of (\ref{conical}).
Thus we can anticipate that
\begin{equation}
( A(\gamma \times \dot{\gamma}), \gamma \times \dot{\gamma} ) = b_1  \dot{\lambda}_1^2 +  b_2  \dot{\lambda}_2^2 \,\,\,.
\end{equation}
In order to find these coefficients $b_1,b_2$, we need   to express
$(\gamma_1,\gamma_2,\gamma_3)$ in terms of $\lambda_1,\lambda_2$. 
This is simply a matter of solving the linear system $$   A (\gamma_1^2,\gamma_1^2,\gamma_1^2)^t = (1,0,0)^t
$$
with
$$
A = \left (\begin {array}{ccc} 1&1&1\\\noalign{\medskip}\left ({\it I_1}-{\it \lambda_1}\right )^{-1}&\left ({\it I_2}-{\it \lambda_1}\right )^{-1}&\left ({\it I_3}-{\it \lambda_1}\right )^{-1}\\\noalign{\medskip}\left ({\it I_1}-{\it \lambda_2}\right )^{-1}&\left ({\it I_2}-{\it \lambda_2}\right )^{-1}&\left ({\it I_3}-{\it \lambda_2}\right )^{-1}\end {array}\right )  \,\,.
$$
Inverting $A$ gives the desired explicit expression,
\begin{equation}  \label{gamaconical}
\left(
 \gamma_1^2 \,,\, \gamma_2^2 \,,\,  \gamma_3^2 
\right ) =
\left (   \frac{(I_1 - \lambda_1)(I_1 - \lambda_2)}{(I_1 - I_3)(I_1 - I_2)}  \,,\, \frac{(I_2 - \lambda_1)(I_2 - \lambda_2)}{(I_2 - I_3)(I_2 - I_1)}  \,,\, \frac{(I_3 - \lambda_1)(I_3 - \lambda_2)}{(I_3 - I_1)(I_3 - I_2)} 
  \right )  \,\,\, .
\end{equation}
It follows that
\begin{equation}
\frac{\partial \gamma}{\partial \lambda_i} = \frac{1}{2} \, (\,\frac{\gamma_1}{\lambda_i - I_1}, \, \frac{\gamma_2}{\lambda_i - I_2}\,,\, \frac{\gamma_3}{\lambda_i - I_3}\,) \,\,\,, \,\,  i = 1,2
\end{equation}
and
\begin{equation}
\gamma \, \times \, \frac{\partial \gamma}{\partial \lambda_i} = \frac{1}{2}\, \left ( \, \frac{\gamma_2 \gamma_3 (I_3-I_2)}{(I_2-\lambda_i)(I_3-\lambda_i)}\,\,,\,\, \frac{\gamma_1 \gamma_3 (I_1-I_3)}{(I_1-\lambda_i)(I_3-\lambda_i)}\,,\, 
\frac{\gamma_1 \gamma_2 (I_2-I_1)}{(I_1-\lambda_i)(I_2-\lambda_i)} \, \right )
\end{equation}
The remaining calculations can be done by brute force, but it was less painful to use a computer algebra system. The results are:
\begin{lemma}
\begin{equation}  \label{geoellipsoid}
(A \dot{\gamma}, \dot{\gamma}) = 1/4\,{\frac {\left ({\it \lambda_2}-{\it \lambda_1}\right ){\it \lambda_1}}{\left ({\it \lambda_1}-{\it I_1}\right )\left ({\it \lambda_1}-{\it I_2}\right )\left ({\it \lambda_1}-{\it I_3}\right )}} \, \dot{\lambda}_1^2 + 1/4\,{\frac {\left ({\it \lambda_1}-{\it \lambda_2}\right ){\it \lambda_2}}{\left ({\it \lambda_2}-{\it I_1}\right )\left ({\it \lambda_2}-{\it I_2}\right )\left ({\it \lambda_2}-{\it I_3}\right )}} \, \dot{\lambda}_2^2
\end{equation}
\begin{equation}  \label{standardmetric}
 || \dot{\gamma}||^2 = 
1/4\,{\frac {{\it \lambda_2}-{\it \lambda_1}}{\left ({\it \lambda_1}-{\it I_1}\right )\left ({\it \lambda_1}-{\it I_2}\right )\left ({\it \lambda_1}-{\it I_3}\right )}} \dot{\lambda}_1^2 +
1/4\,{\frac {{\it \lambda_1}-{\it \lambda_2}}{\left ({\it \lambda_2}-{\it I_1}\right )\left ({\it \lambda_2}-{\it I_2}\right )\left ({\it \lambda_2}-{\it I_3}\right )}} \, \dot{\lambda}_2^2 
\end{equation}
\begin{equation} \label{twisted}
( A(\gamma \times \dot{\gamma}), \gamma \times \dot{\gamma} ) = b_1  \dot{\lambda}_1^2 +  b_2  \dot{\lambda}_2^2 = \frac{1}{4} \frac{\lambda_2 (\lambda_2 - \lambda_1)}{(\lambda_1 -I_1)(\lambda_1 -I_2)(\lambda_1 -I_3)}\,\dot{\lambda}_1^2
+ \frac{1}{4} \frac{\lambda_1 (\lambda_1 - \lambda_2)}{(\lambda_2 -I_1)(\lambda_2 -I_2)(\lambda_2 -I_3)} \,  \dot{\lambda}_2^2 \,\,\,.
\end{equation}
\end{lemma}
Equation (\ref{geoellipsoid}) indeed tell us that the ellipsoid geodesics problem separates
(the standard metric $\dot{x}^2$ in the ellipsoid $ \sum x_j^2/I_j = 1$
is equivalent to the deformed metric $(A \dot{\gamma}, \dot{\gamma})$ in the sphere
$ \sum \dot{\gamma}_j^2 = 1\,\,,\,\, \gamma_j = x_j/\sqrt{I_j} \,$).
The standard metric (\ref{standardmetric}) on the sphere, which gives a superintegrable system,
also appears separable in conical coordinates (but in a somewhat cumbersome format). 



 Let us now  write the symplectic form (\ref{symplecticform}) using
$\lambda_1,\lambda_2,\dot{\lambda_1},\dot{\lambda_2}$ as coordinates.
This requires expressing the following objects in terms of those:
\begin{lemma}
$$ dS^2 = (\gamma, \frac{\partial \gamma}{\partial \lambda_1} \times \frac{\partial \gamma}{\partial \lambda_2} ) \, d\lambda_1 \wedge d\lambda_2 
= \frac{\lambda_2-\lambda_1}{\sqrt{(\lambda_1-I_1)(I_2-\lambda_1)(I_2-\lambda_1)(\lambda_2-I_1)(\lambda_2-I_2)(I_3-\lambda_2)}}  d\lambda_1 d\lambda_2 \,\,.
$$
\begin{eqnarray*} (\gamma, A(\dot{\gamma} \times \gamma)) & = & \frac{1}{2} \, \gamma_1 \gamma_2 \gamma_3 \left( \dot{\lambda}_1 ( \frac{I_1(I_2-I_3)}{(I_2-\lambda_1)(I_3-\lambda_1)} + {\rm cyclic} ) +
\dot{\lambda}_2 ( \frac{I_1(I_2-I_3)}{(I_2-\lambda_2)(I_3-\lambda_2)} + {\rm cyclic} ) \right)  \\
& = & \frac{1}{2} \sqrt{\frac{(\lambda_2-I_1)(\lambda_2-I_2)(I_3-\lambda_2)}{(\lambda_1-I_1)(I_2 -\lambda_1)(I_3-\lambda_1)}} \dot{\lambda}_1 + \frac{1}{2} \sqrt{\frac{(\lambda_1-I_1)(I_2-\lambda_1)(I_3-\lambda_1)}{(\lambda_2-I_1)(\lambda_2-I_2)(I_3-\lambda_2)}} \dot{\lambda}_2
\end{eqnarray*}
 $$
\frac{\gamma_2^2}{I_2} + \frac{\gamma_3^2}{I_3} = \frac{\lambda_1 \lambda_2}{I_1 I_2 I_3}
$$
$$
 \frac{\gamma_2^2 + \gamma_3^2}{I_2I_3} + \frac{\gamma_1^2 + \gamma_3^2}{I_1I_3} + \frac{\gamma_1^2 + \gamma_2^2}{I_1I_2} = \frac{\lambda_1 + \lambda_2}{I_1I_2I_3}
$$
$$  (J,K) = \frac{1}{2} (1-\frac{b}{a}) (1+\frac{b}{a})^2 (\lambda_2-\lambda_1)
\left( \frac{\dot{\lambda}_1}{(\lambda_1-I_1)(\lambda_1-I_2)(\lambda_1-I_3)} +
\frac{\dot{\lambda}_2}{(\lambda_2-I_1)(\lambda_2-I_2)(I_3-\lambda_2)} 
\right)
$$
\begin{equation}
F(\gamma) = (1+\frac{b}{a})^2 \, (\lambda_1 + \mu b^2) (\lambda_2 + \mu b^2)
\end{equation}
\begin{equation}
  T = \frac{1}{8}  (1+\frac{b}{a})^2 (\lambda_2-\lambda_1) \left(
\frac{\lambda_2 + \mu b^2}{(\lambda_1-I_1)(I_2-\lambda_1)(I_3-\lambda_1)} \dot{\lambda}_1^2  + 
\frac{\lambda_1 + \mu b^2}{(\lambda_2-I_1)(\lambda_2-I_2)(I_3 - \lambda_2)} \dot{\lambda}_2^2 
\right)
\end{equation}
\begin{equation} p_1 = \partial T/\partial \dot{\lambda}_1 = \frac{1}{4}  (1+\frac{b}{a})^2 (\lambda_2-\lambda_1) 
\frac{\lambda_2 + \mu b^2}{(\lambda_1-I_1)(I_2-\lambda_1)(I_3-\lambda_1)} \dot{\lambda}_1 
\end{equation}
\begin{equation}
p_2 =  \partial T/\partial \dot{\lambda}_2 =  \frac{1}{4}  (1+\frac{b}{a})^2 (\lambda_2-\lambda_1) \frac{\lambda_1 + \mu b^2}{(\lambda_2-I_1)(\lambda_2-I_2)(I_3 - \lambda_2)} \dot{\lambda}_2
\end{equation}
\end{lemma}
Putting everything together, we get:
\begin{proposition}  Darboux coordinates in new time:
\begin{eqnarray} \omega_{conformal} & = &[(\lambda_1+\mu b^2)(\lambda_2+\mu b^2)]^{\frac{b-a}{2a}} \omega_{nh} = \\
& = & [(\lambda_1+\mu b^2)(\lambda_2+\mu b^2)]^{\frac{b-a}{2a}} \left( dp_1 d\lambda_1 + dp_2 d\lambda_2 + \frac{1}{2}(1-b/a) (\frac{p_1}{\lambda_2+\mu b^2} - \frac{p_2}{\lambda_1+\mu b^2}) d\lambda_1 d\lambda_2 \right) \nonumber \\ 
 & = & d(P_1 d\lambda_1 + P_2 d\lambda_2) \,\,\,\,\,\, \,\,\,\, {\rm  Darboux!}
 \end{eqnarray}
\end{proposition}
\begin{proposition}  The Hamiltonian in spheroconical coordinates: in new time 
\begin{equation} d\tau/dt = 
F(\gamma)^{\frac{b-a}{2a}} \,\,\,,\,\,\, F(\gamma) = (1+\frac{b}{a})^2 \, (\lambda_1 + \mu b^2) (\lambda_2 + \mu b^2)  
\end{equation}
defining
\begin{equation}
P_1 = [(\lambda_1+\mu b^2)(\lambda_2+\mu b^2)]^{\frac{b-a}{2a}} p_1  \,\,,\,\, P_2 = [(\lambda_1+\mu b^2)(\lambda_2+\mu b^2)]^{\frac{b-a}{2a}} p_2  \,\,.
\end{equation}
the Hamiltonian  becomes
\begin{equation} 2H = \frac{p_1^2}{c_1} +  \frac{p_2^2}{c_2} =  \frac{P_1^2}{[(\lambda_1+\mu b^2)(\lambda_2+\mu b^2)]^{\frac{b-a}{a}} c_1} +  \frac{P_2^2}{[(\lambda_1+\mu b^2)(\lambda_2+\mu b^2)]^{\frac{b-a}{a}}c_2} 
\end{equation}
with
\begin{equation}
c_1 = \frac{1}{4}(1+ \frac{b}{a})^2 (\lambda_2-\lambda_1)  \frac{\lambda_2 + \mu b^2}{(\lambda_1-I_1)(I_2-\lambda_1)(I_3-\lambda_1)}\,\,,\,\, c_2 = \frac{1}{4}(1+ \frac{b}{a})^2 (\lambda_2-\lambda_1)  \frac{\lambda_1 + \mu b^2}{(\lambda_2-I_1)(\lambda_2-I_2)(I_3-\lambda_2)} \,\,.
\end{equation}
\end{proposition}


\section{Some  directions for research}  

\subsection{Nonintegrability.}

If the rolling ball has  its mass concentrated at the geometric center,
then the inertias $I_j$ vanish. Hence,  the problem with $I_j = \epsilon \, \bar{I}_j $, is a perturbation of
the geodesics on the sphere.  As a possible attempt to prove nonintegrability rigorously 
is to compute the monodromy
of the linearized equations along a great circle (solution of the unperturbed problem) and
 apply Melnikov or Moralis-Ramis methods. Parhaps this may be done directly
in $\gamma \in \R^3$.  In   sphero-conical coordinates,
using rescaled variables $ \lambda_j = \epsilon \bar{\lambda}_j$
and new time with factor
\begin{equation}
F(\gamma) = (1+\frac{b}{a})^2 \, (\epsilon \lambda_1 + \mu b^2) (\epsilon \lambda_2 + \mu b^2)
\end{equation}
we obtain  (dropping the bars) the following Hamiltonian perturbation
of the geodesics on the sphere:
\begin{equation} 2H =   \frac{P_1^2}{[(\epsilon \lambda_1+\mu b^2)(\epsilon \lambda_2+\mu b^2)]^{\frac{b-a}{a}} a_1} +  \frac{P_2^2}{[(\epsilon \lambda_1+\mu b^2)(\epsilon \lambda_2+\mu b^2)]^{\frac{b-a}{a}}a_2} 
\end{equation}
with
\begin{equation}
a_1 = \frac{1}{4}(1+ \frac{b}{a})^2 (\lambda_2-\lambda_1)  \frac{\epsilon \lambda_2 + \mu b^2}{(\lambda_1-I_1)(I_2-\lambda_1)(I_3-\lambda_1)}\,\,,\,\, a_2 = \frac{1}{4}(1+ \frac{b}{a})^2 (\lambda_2-\lambda_1)  \frac{\epsilon \lambda_1 + \mu b^2}{(\lambda_2-I_1)(\lambda_2-I_2)(I_3-\lambda_2)}  \,\,.
\end{equation}

\subsection{Rolling with sliping but no twisting and Rattleback phenomena}

From  the mathematical point of view, it is   legitimate to study systems subject to the no-twist constraint,  sliding allowed. 
The question is whether these constraints could be (even approximately) realized in practice.  Nonholonomic constraints are approximately realized as limits
of certain collisional systems. It would be interesting to conceive
a ``peg-leg'' version for no-sliping and the no-twist constraints\footnote{We thank Andy Ruina for a discussion on this subject; see 
$ http://ruina.tam.cornell.edu/ $. }.

 We believe that the dynamical torque (\ref{dyntorque}) the heart of the rattleback rotation reversal phenomenon.  
There is a extensive literature on the topic, but we would like to see
an explanation based on the joint geometries of
$\Sigma_2$, say with $\Sigma_1$ being the plane. 
Is there a rubber rattleback?  Here the idea is not looking at rotations (which are forbidden) but at translation reversal.

\subsection{Skiding dynamics: mathematical hockey}

  As far as we know holonomic systems  $ T(q_1,R,\dot{q}_1,\dot{R})$  
  given by (\ref{lagrangianofQ}),  describing the
motion of a body ${\cal B}$  {\it skiding}  over a surface $\Sigma_1$ have  not been studied systematically.    
  In general the rotational and surface dynamics should be  coupled and some interesting effects of mathematical hockey may be 
 anticipated\footnote{Hockey (see  http://www.exploratorium.edu/hockey/)  could be even more exciting with a strictly convex puck.}. 
A first project in this
direction is to explore symmetries (as we  discussed briefly).

\subsection{Rubber rolling of surfaces of revolution  over a sphere }

  Conditions under which a generalized Chaplygin system has an invariant measure
are given in \cite{Cantrijn} .
It would be interesting  to  make explicit which condition results for marble
rolling of an arbitrary surface $\Sigma_2$ over a sphere or a plane.    Also,
implementating proposition \ref{reductionchaplygin}  on  exemples such as bodies of revolution and ellipsoids is in order,
as well as the higher dimensional extensions \cite{FedorovKozlov,Jovanovic}.  
These situations have been  studied in the marble rolling context,
see tables 1 and  2  of \cite{BM}.
On another talk, rubber rolling of an homogenous sphere over arbitrary surfaces $\Sigma_1$, see \cite{BMK} for the
marble counterpart.

\subsection{Suslov rolling}

While the usual no-slip constraint is enforced, the no-twist constraint $\omega_3 = 0 $ is replaced by $\Omega_3 = 0$,
that is,  one imposes that the third component of angular velocity {\it in the body frame} should vanish.  This constraint can be implemented by a inner spherical support.

\subsection{Rubber rolling and Yang-Mills }

Among all convex bodies ${\cal B}$ with  the same area, rolling over a sphere,  which one  minimizes the
Yang-Mills functional? This functional, important for gauge theories, is given by
 $$ \int_{S^2} \,  ||K||^2  \, dA  \,\,\, , $$ where $K$ is the curvature of the connection.   Here we restric ourselves to  rubber connections $\phi_{{\cal B}}$
  on the bundle  $ SO(3) \times S^2 $. A nice exercise is to  compute the rubber connection 1-form and its curvature, for a
tri-axial ellipsoid rolling over a sphere.

A reasonable conjecture  is that the sphere minimizes the YM functional.  Clearly, the
minimum attains the value zero when the area of {\cal B} is the same as
the area of $S_a$.  

   Yang Mills connections
over Riemann surfaces  were considered on an important paper by Atiyah and  Bott \cite{AtiyahBott}. 
Particularize to rolling problems when the base is a surface of constant negative curvature.  


\bigskip
\bigskip

\noindent {\bf Acknowledgements.}   Both authors thank
 Prof. Mikhail A. Sokolovskiy for inviting us to the IUTAM symposium 2006 on
 Hamiltonian Dynamics, Vortex Structures, Turbulence, at the Steklov Institute, Moscow, 
 $http://conf2006.rcd.ru/$ where me met Alexei Borisov, Ivan Mamaev, 
 Alexander Kilin, Yuri Fedorov, and renewed contact with Anatoly Neishtadt.
 This meeting was a wonderful experience. JK slso thanks  Richard Murray and Jerry Marsden
for hosting a Fulbright visit, Winter 2005, where this research started: a great time in company of  
 Tudor Ratiu, Marco Castrillon, and Hernan Cendra.
 


\begin{thebibliography}{99}  


\bibitem{Arnold}
Arnold, V.~I., {\em Mathematical Methods of Classical Mechanics},
volume~60 of {\em Graduate Texts in Math.}
\newblock Springer-Verlag, First Edition 1978, Second Edition, 1989.

\bibitem{Arnoldetal}
Arnold, V.~I.,  Kozlov, V.V.,  A.~I. Neishtadt, A.I.,
{\em Mathematical aspects of classical and celestial mechanics},
in Arnold, V.~I., editor, {\em Dynamical Systems III, Encyclopaedia of Mat. Sciences},
\newblock Springer-Verlag, 1988.

\bibitem{AtiyahBott} Atiyah, M., Bott, R., 
The Yang-Mills Equations over Riemann Surfaces, 
{\em Philosophical Transactions of the Royal Society of London. Series A, Mathematical and Physical Sciences}, \textbf{308}:1505, 1983.
 



\bibitem{Blackall}  
Blackall, C.J., On volume integral invariants of non-holonomic dynamical systems,
{\em Amer. J. Math.} \textbf{63}:1, 155--168, 1941.

\bibitem{Bloch}
Bloch, A.M., {\em  Nonholonomic mechanics and control},
\newblock  Springer-Verlag, 2003.

\bibitem{BKMM}
Bloch, A.M., Krishnaprasad, P.S., Marsden, J.E., Murray, R.M., Nonholonomic mechanical
systems with symmetry, {\em Arch. Ratl. Mech. Anal.}, \textbf{136}, 21--99.1996.

\bibitem{MontBor} Bor, G., and Montgomery, R., $G_2$ and the ``rolling distribution'', preprint.

\bibitem{Borisovsphere}
Borisov, A. V., Mamaev, I.S., Chaplygin's ball rolling problem is Hamiltonian,
{\em Mathematical Notes (Matematicheskie Zametki)}, \textbf{70}:5, 793-795, 2001.

\bibitem{Borisovsphere1}
Borisov, A. V., Mamaev, I.S., Obstacle to the Reduction of Nonholonomic Systems
to the Hamiltonian Form, {\em Doklady Physics USSR}, \textbf{47}:12, 892--894, 2002.

\bibitem{BM}
Borisov, A. V., Mamaev, I. S.,
The Rolling Body Motion Of a Rigid Body on a Plane and a 
Sphere. Hierarchy of Dynamics,
{\em Regular and Chaotic Dyn.}, \textbf{7}:2, 177--200, 2002.

\bibitem{BMK}
Borisov A.V., Mamaev I.S. and Kilin A.A., The rolling motion of a ball on
a surface. New integrals and hierarchy of dynamics, {\em Regular and Chaotic 
Dyn.}, \textbf{7}:2, 201--218, 2002.

\bibitem{BM1}  Borisov, A.V., Mamaev, I.S.,  Isomorphism and Hamilton representation of some nonholonomic systems, {em Siberian Mathematical Journal}, \textbf{48}:11,2007, to appear.
see preliminary version  at arXiv: nlin.SI/0509036 v1 (Sept. 21 2005).

\bibitem{personalcommunication} Borisov, A.V, Mamev, I.S.,  this issue.

\bibitem{Bryant} Bryant, R., Hsu,L., Rigidity of integral curves of rank 2 distributions, Inventiones Mathematicae \textbf{114}, 435--461,  1993.

\bibitem{Cantrijn1}
Cantrijn, F., de L\'eon, M., Marrero, J.C., de Diego, D., 
Reduction of nonholonomic mechanical systems with symmetries,
{\em Rep. Math. Phys.}  \textbf{42}:1/2, 25--45, 1998.

\bibitem{Cantrijn3}
Cantrijn, F.,  de L\'eon, M., de Diego, D., On almost-Poisson structures in nonholonomic
mechanics, {\em Nonlinearity}  \textbf{12}, 721--737, 1999.

\bibitem{Cantrijn}
Cantrijn, F., Cort\'es, J., de L\'eon, M., de Diego, D., On the geometry of generalized
Chaplygin systems, {\em Math. Proc. Camb. Phil. Soc.}  \textbf{132}, 323--351, 2002.


\bibitem{Cartan}  
Cartan, \'E. , Le syst\`emes de Pfaff \`a cinq variables et les \'equations aux d\`eriv\`ees partielles du second ordre, 
 {\em Ann. Sci. \`Ecole Normale } \textbf{27}:(3), 109--192, 1910.


\bibitem{Cartan1} 
Cartan, \'E., Sur la repres\'entation
g\'eom\'etrique
des syst\`emes mat\'eriels non holonomes, {\em Proc. Int. Congr. Math. Bologna} \textbf{4},    253--261, 1928.


\bibitem{CeMaRa}
Cendra, H., Marsden, J. E., Ratiu,  T. S., Geometric mechanics, Lagrangian reduction and nonholonomic systems, 
{\em Mathematics Unlimited-2001 and Beyond}, (B. Enguist and W. Schmid, eds.), Springer-Verlag, New York, 221--273, 2001.

\bibitem{Chaplygin0}
Chaplygin, S.A., On the theory of the motion of nonholonomic systems. Theorem
on the reducing factor, {\em Mat. Sbornik}, \textbf{28}, 303--314, 1911. 


\bibitem{Chaplygin1}
Chaplygin, S.A., On a ball's rolling on a horizontal plane, {\em Reg. Chaot. Dyn.},
\textbf{7}:2, 131--148, 2002. \\  Original paper in {\em Math. Sbornik} \textbf{24}, 139--168, 1903.

\bibitem{Cortes1}
Cort\'es, J., {\em Geometric, Control and Numerical Aspects of Nonholonomic Systems},
 \newblock Springer-Verlag, 2002. 

\bibitem{CortesManolo}
Cort\'es, J., de L\'eon, M., de Diego, D., Mart\'{\i}nez, S., Geometric description of vakonomic and nonholonomic dynamics; comparison of solutions,
{\em SIAM J. Control Optim.} \textbf{41}:5, 1389--1412, 2003.

\bibitem{Cushmanbook}
Cushman, R., Bates, L., {\em Global aspects of Classical Integrable Systems},
\newblock Birkh\"auser, Basel, 1997.



\bibitem{Duis} 
Duistermaat, J.J. [2000], Chapygin's sphere,   in R. Cushman, J. J. Duistermaat and J. \'Sniatycki: {\em Chaplygin and the Geometry of
Nonholonomically Constrained Systems}, 2000 (in preparation).  See arxiv:   


\bibitem{Kurt2}  Ehlers, K.,  Koiller, J., Montgomery, R., Rios, P.M.,  Nonholonomic systems via moving frames: Cartan equivalence
and Chaplygin hamiltonization, In The Breadth of Symplectic and Poisson Geometry, Festschrift in Honor of Alan Weinstein, ed. by
J.Marsden and T.Ratiu, Birkhäuser, Boston, 2005

\bibitem{EhlersIUTAM} Ehlers, K., Koiller, J., Rubber rolling: geometry and dynamics of 2-3-5 distributions, in Proceedings IUTAM symposium 2006 on
 Hamiltonian Dynamics, Vortex Structures, Turbulence, at the Steklov Institute,
 Moscow,
to appear.

\bibitem{EhlersKoiller}  Ehlers, K., Koiller, J., Classification of the 2-3-5 nonholonomic geometries, in preparation.


\bibitem{FedorovKozlov} 
Fedorov, Yu. N., Kozlov, V.V., Various aspects of n-dimensional rigid body dynamics, in Kozlov, V.V. (editor) {\em Dynamical Systems in
Classical Mechanics}, volume 168 of {\em AMS Translations series 2}, 1995.

\bibitem{Fedorov}
Fedorov, Yu. N.,  Jovanovic, B.,
Nonholonomic LR systems as Generalized Chaplygin systems with an Invariant Measure and Geodesic Flows on Homogeneous Spaces, 
J. Nonlinear Science, \textbf{14}:1, 341--381, 2004 ( arxiv.org/abs/math-ph/0307016, 2003).

\bibitem{Haller} 
Haller, S., Rybicki, T.,  On the group of diffeomorphisms preserving a locally conformal symplectic structure,
 {\em Ann. Global Anal. Geom.} \textbf{17}, 475--502, 1999.

\bibitem{Haller1}  
Haller, S., Rybicki, T., Symplectic reduction for locally conformal symplectic manifolds, {\em J. Geom. Phys.} , \textbf{37}, 262--271, 2001. 



\bibitem{Hertz}
Hertz, H., {\em
The principles of mechanics presented in a new form by Heinrich Hertz, with an introduction by H. von Helmholtz}, 
\newblock Macmillan,  London, New York,  1899. 


\bibitem{Iliyev}  
Iliyev, Il., On the conditions for the existence of the reducing Chaplygin factor,
{\em P.M.M. USSR} \textbf{49}:2, 295--301, 1985.

\bibitem{Jovanovic} 
Jovanovic, B. , Some multidimensional integrable cases of nonholonomic
rigid body dynamics, {\em Regular $\&$ Chaotic Dynamics} \textbf{8}:1, 125--132, 2003.

\bibitem{Kilin} Kilin, A.A.,  The Dynamics of Chaplygin ball: the qualitative and computer analysis,{\em  Reg. Chaotic Dynamics} 
\textbf{ 6}:3, 291--306, 2001.  

\bibitem{KobayashiOliva}
Kobayashi, M.H., Oliva, W.M., A note on the conservation of energy and volume in the
setting of nonholonomic mechanical systems, {\em Qual. Theory Dyn. Systems}, \textbf{5}, 255-283, 2004. 

\bibitem{KoillerArma}
Koiller, J., Reduction of some classical non-holonomic systems with symmetry, {\em Arch. Rational Mech. Anal.}  \textbf{118},  113--148, 1992.


\bibitem{KoillerRomp}  
Koiller, J., Rios, P.M., Ehlers, K., Moving frames for cotangent bundles, 
{\em Rep.   Math.  Phys.},  \textbf{49}:2/3, 225--238, 2002.



\bibitem{Kozlov}
Kozlov, V.V., On the integration theory of equations of nonholonomic mechanics,
{\em Reg. Chaot. Dyn.}, \textbf{7}:2, 161--176, 2002. 


 
\bibitem{KM}
Koon, W. S., Marsden, J.E., The Hamiltonian and Lagrangian approaches to the
dynamics of nonholonomic systems, {\em Rep. Math. Phys.} , \textbf{40}, 21--62, 1997.

\bibitem{KM1}
Koon, W. S., Marsden, J.E., Poisson reduction for nonholonomic mechanical systems with symmetry, 
{\em Reports on Math. Phys.}, \textbf{42}, 101--134, 1998. 

\bibitem{KupkaOliva}
Kupka, I., Oliva, W.M., The Non-Holonomic Mechanics, {\em J. Diff. Equations} \textbf{169}, 169--189, 2001.

\bibitem{Levi}
Levi, M.[1996], Composition of rotations and parallel transport, {\em Nonlinearity} \textbf{9}, 413--419, 1996.

\bibitem{Marle1}
Marle, C.M., Various approaches to conservative and nonconservative nonholonomic systems,  {\em Rep. Math. Phys.} \textbf{42}:1/2, 211--229,
1998.

\bibitem{Marle2}
Marle, C.M., On symmetries and constants of motion in Hamiltonian systems with nonholonomic constraints, 
{\em Banach Center Publ.}, \textbf{59}, 223--242, 2003.


\bibitem{Marle3} Marle, C.M., From momentum maps and dual pairs to symplectic and Poisson groupoids,
 In: {\em The Breadth of Symplectic and Poisson Geometry}, Festschrift in Honor of Alan Weinstein, Birkhauser, Boston, ed. by J. Marsden and T. Ratiu,  493-523, 2005.
 
 

\bibitem{MR} Marsden, J.E., and Ratiu, T., {\em Introduction to Mechanics and Symmetry, 
Texts in Applied Mathematics}  \textbf{vol. 17}, 2nd ed. , Springer, New York, 1999.

\bibitem{MW}
Marsden, J.E., Weinstein, A.,  Reduction of symplectic manifolds with symmetry,
{\em Rep. Math. Phys.} \textbf{5}, 121--130, 1974.

\bibitem{Montgomery} Montgomery, R., {\em A Tour of Subriemannian Geometries, their Geodesics, and Applications, AMS Surveys and
Monographs}  \textbf{91}, 2002.

\bibitem{Naranjo}  Naranjo, L.G., Ph.D. thesis, University of Arizona
(see $http://math.arizona.edu/~luisg/ $).

\bibitem{NF} 
Neimark, J.I.,  N.~A.~Fufaev, N.A., {\em Dynamics of nonholonomic systems}
\newblock volume~33 of {\em AMS Translations of Mathematical Monographs}, Providence, 1972.

\bibitem{Oliva}
Oliva, W.M., {\em Geometric Mechanics}, volume 1798 of {\em Springer Lecture Notes
in Mathematics}
\newblock Springer Verlag, 2002.





 
 \bibitem{Schaft} van der Schaft, A.J., Port-Hamiltonian systems: an approach to modeling and control of complex physical systems, Proceedings of the Sixteenth International Symposium on Mathematical Theory of Networks and Systems (MTNS2004), Leuven, Belgium, July 5--9, 2004.

\bibitem{Sniatycki98}
\'Sniatycli, J., Nonholonomic Noether theorem and reduction of symmetries,
{\em Rep. Math. Phys.} \textbf{42}:1/2, 5--23, 1998.




\bibitem{Stanchenko}
Stanchenko, S.V., Non-holonomic Chaplygin systems, {\em P.M.M. USSR} \textbf{53}:1,
11--17, 1985.

\bibitem{Struik} Struik, D. J., Lectures on classical differential geometry, Addison-Wesley, 1950.


\bibitem{Vaisman} 
Vaisman, I., Locally conformal symplectic manifolds. 
{\em Internat. J. Math. Math. Sci.} \textbf{8}:3, 521--536, 1985. 

\bibitem{Veselovs1}
Veselov, A. P., Veselova, L. E.,  Flows on Lie groups with a nonholonomic constraint and integrable non-Hamiltonian systems. (Russian), {\em Funktsional. Anal. i Prilozhen.} \textbf{20}:4, 65--66. 
English translation: {\em Functional Anal. Appl.} \textbf{20}:4,308--309, 1986.

\bibitem{Veselovs}
Veselov, A. P., Veselova, L. E., Integrable nonholonomic systems on Lie groups. (Russian) {\em Mat. Zametki}
 \textbf{44}:5   604--619, 701; translation in {\em Math. Notes} \textbf{44}:5/6, 810--819, 1989. 
 
 
\bibitem{Yoshimura} Yoshimura, H., Marsden, J.,
Dirac Structures and Implicit Lagrangian Systems in Electric Networks,
17th International Symposium on
Mathematical Theory of Networks and Systems,   Kyoto,   2006.

\bibitem{AissaWade} 
Wade, A., Conformal Dirac structures, {\em Lett. Math. Phys.}
\textbf{53},  331--348, 2000.

\bibitem{Blochnonlinearity}  
Zenkov, D.V., Bloch, A.M.,   Invariant measures of nonholonomic flows with
internal degrees of freedom, {\em Nonlinearity} \textbf{16}, 1793--1807, 2003.





\end{thebibliography}
\end{document}